\documentstyle[10pt]{article}
\begin{document}

\baselineskip 15pt
\parindent=1em
\hsize=12.3 cm \textwidth=12.3 cm
\vsize=18.5 cm \textwidth=18.5 cm

\def\supt{{\rm supt}}
\def\dom{{\rm dom}}
\def\bfone{{\bf 1}}
\def\Gen{{\rm Gen}}
\def\rk{{\rm rk}}
\def\cali{{\cal I}}
\def\calj{{\cal J}}
\def\Spec{{\rm Spec}}

\title{How special is your Aronszajn tree?}
\author{
Chaz Schlindwein \\
Division of Mathematics and Computer Science \\
Lander University \\
Greenwood, South Carolina 29649, USA\\
{\tt chaz@lander.edu}}

\maketitle

\def\forces{\mathbin{\parallel\mkern-9mu-}}
\def\notforces{\,\nobreak\not\nobreak\!\nobreak\forces}

\def\restr{\,\hbox{\vrule height8pt width.4pt depth0pt
   \vrule height7.75pt width0.3pt depth-7.5pt\hskip-.2pt
   \vrule height7.5pt width0.3pt depth-7.25pt\hskip-.2pt
   \vrule height7.25pt width0.3pt depth-7pt\hskip-.2pt
   \vrule height7pt width0.3pt depth-6.75pt\hskip-.2pt
   \vrule height6.75pt width0.3pt depth-6.5pt\hskip-.2pt
   \vrule height6.5pt width0.3pt depth-6.25pt\hskip-.2pt
   \vrule height6.25pt width0.3pt depth-6pt\hskip-.2pt
   \vrule height6pt width0.3pt depth-5.75pt\hskip-.2pt
   \vrule height5.75pt width0.3pt depth-5.5pt\hskip-.2pt
   \vrule height5.5pt width0.3pt depth-5.25pt}\,}

   \centerline{{\bf Abstract}}

   We answer a question of Shelah by constructing a model of
    Suslin's Hypothesis in which
   there is an Aronszajn tree $T$ such that for every unbounded
    $E\subseteq\omega_1$ we have that
   $T$ is not $E$-special.  We may require that CH holds, or
    that CH fails, or that Kurepa's hypothesis holds or fails, or that
   there is a stationary $S\subseteq\omega_1$ such that every
    Aronszajn tree is
   $S$-$*$-special, or other variants.

\eject

\section{Introduction}

In this note, we answer the following questions of Shelah
[Sh, Remark IX.4.9(5)]:

\proclaim Q1. Is it consistent that Suslin's Hypothesis holds yet
there is an
 Aronszajn tree $T$ such that for every unbounded $E\subseteq\omega_1$
 we have
that $T$ is not $E$-special?

\proclaim Q2. Is it consistent that there is a stationary $S\subseteq\omega_1$
 such that every Aronszajn tree is $S$-$*$-special yet there is an
  Aronszajn tree $T$ such that
for every unbounded $E\subseteq\omega_1$ we have that $T$ is not $E$-special?

Notice that a positive answer to Q2 yields a positive
answer to Q1 (see [PIF, Claim IX.3.4(1)] or Lemma 2.9, below).

The intent of Shelah's questions is to delineate the distinction
 between two different
notions of ``special'' for Aronszajn trees.  Shelah [Sh, Remark IX.4.9(2)]
addresses this distinction as follows:

\proclaim  Theorem (Shelah).
If ZFC is consistent, then so is ZFC plus there is an unbounded
$E\subseteq\omega_1$
such that every Aronszajn
tree  is $E$-special and there is an Aronszajn tree $T$ such that for every
stationary $S\subseteq\omega_1$ we have that $T$ is not $S$-$*$-special.

This shows in a strong way
that $E$-specialness does not entail $S$-$*$-specialness.  Thus it
is natural to consider whether $S$-$*$-specialness entails $E$-specialness.
This is the motivation for the two questions of Shelah given above.
We show that
the answer to each of the above questions is positive; that is,
neither version of ``specialness'' implies the other.  
In particular, we demonstrate the following theorem.

\proclaim Theorem.  If ZFC is consistent, then so is ZFC 
plus there is a stationary $S\subseteq\omega_1$
such that every Aronszajn tree is $S$-$*$-special and there
is an Aronszajn tree $T$
such that for every unbounded $E\subseteq\omega_1$ we have
that $T$ is not $E$-special.

We produce models with certain additional properties.  For example,
we give a model exemplifying a positive answer to Q1
in which $T$ has no stationary antichain.
The models we use are variations of
models that have appeared in [S], [S], [S].
We show that each of these models satisfies the statement: ``there is no
unbounded $E\subseteq\omega_1$ such that $T$ is $E$-special (in some of the
cited papers, the distinguished Aronszajn tree is denoted $T^*$ rather than $T$).''  The
innovation which leads to our answering Shelah's questions
 is our formulation
of a new preservation property (see Definition 3.6).  Any forcing that
satisfies this property maintains the non-$E$-specialness
 (for every unbounded $E\subseteq\omega_1$) of some
 appropriate Suslin tree
  of the ground model.

\proclaim Notations.  We say that a map $f$ is order-preserving iff $x\leq y$ implies $f(x)\leq f(y)$, whereas we say that $f$ is
strictly order-preserving iff $x<y$ implies $f(x)<f(y)$.  For $T$ a tree we let $T_\beta=\{x\in T\,\colon\allowbreak
{\rm rk}(x)=\beta\}$ and for $E$ a set of ordinals we let 
$T_E=\bigcup\{T_\beta\,\colon\allowbreak\beta\in E\}$ (the fact that
 each ordinal is
literally a set of ordinals renders this ambiguous, but it is
 always clear in
context).
 We say that $T$ is an $\omega_1$-tree iff every level of\/ $T$ is countable and for every $x\in T$, whenever ${\rm rk}(x)<\beta<\omega_1$ then there are at least two successors of $x$ in $T_\beta$, and each node whose rank is a limit ordinal is uniquely determined by its set of predecessors.

\section{Specializations of $\omega_1$-trees}

There are various ways in which an Aronszajn tree $T$ can be specialized.  
The classical notion is that $T$ is special iff there is a strictly order-preserving 
function from $T$ into ${\bf Q}$.  The essential point is that a special Aronszajn tree 
cannot be a Suslin tree.  Baumgartner [B] and Shelah [S, chapter IX] investigate weaker 
notions of ``special'' that also ensure non-Suslinity.  Of particular interest are the 
following two definitions.

\proclaim Definition 2.1.  Suppose $T$ is an Aronszajn tree and 
$E\subseteq\omega_1$ is unbounded.  We say that $T$ is $E$-special iff 
there is a strictly order-preserving map from $T_E$ into ${\bf Q}$.

\proclaim Definition 2.2.  Suppose $T$ is an Aronszajn tree and $S$ is a  subset of $\omega_1$ consisting of limit ordinals.  We say that 
$T$ is $S$-$*$-special iff there is a function $f$ mapping $T_S$ into $\omega_1$ such that $(\forall x\in T_S)\allowbreak(f(x)<{\rm rk}(x))$
and whenever $x<y$ are in $T_S$ then $f(x)\ne f(y)$.

\proclaim Lemma 2.3.  Suppose $T$ is an Aronszajn tree and either $T$ is 
$E$-special for some unbounded $E$ or $T$ is $S$-$*$-special for some 
 stationary $S$ consisting of limit ordinals.  Then $T$ is not Suslin.

Proof:  In the first case, let $f$ from $T_E$ into ${\bf Q}$ be a specializing function.  
For some $r\in {\bf Q}$ we have that $f^{-1}(r)$ is uncountable; necessarily 
$f^{-1}(r)$ is an antichain of $T$.  In the second case, suppose 
$f$ is as in Definition 2.2. For each $\alpha\in S$ let $g(\alpha)=
{\rm min}\{f(x)\,\colon\allowbreak x\in T_\alpha\}$.  By Fodor's theorem we may choose 
$\gamma$ such that $g^{-1}(\gamma)$ is uncountable.  Clearly $\{x\in T_S\,\colon\allowbreak
f(x)=\gamma\}$ is an uncountable antichain of $T$.

In the next two Definitiona the two notions of ``special'' introduced above
are extended to
$\omega_1$-trees  (see also [Sch]).  We show in Lemma 2.9 that
an $\omega_1$-tree that is special in either of these two extended senses is neither 
Suslin nor Kurepa.
We also show in Lemma 2.8 that for Aronszajn trees,
the two extended definitions coincide with the earlier definitions.  Until then, 
we shall specify, e.g., ``$E$-special in the sense of Definition 2.4.''

\proclaim Definition 2.4.  Suppose\/ $T$ is an\/ $\omega_1$-tree and\/ $E\subseteq
\omega_1$ is uncountable.  We say that\/ $T$ is $E$-special iff there is an
order-preserving\/ $f$ mapping\/ $T_E$ into\/ ${\bf Q}$ such that whenever\/
$\{x,y,z\}\subseteq T_E$ and\/ $f(x)=f(y)=f(z)$ and\/ $x<y$ and\/ $x<z$, 
then\/ $y$ and\/ $z$ are comparable.

\proclaim Definition 2.5. Suppose\/ $T$ is an\/ $\omega_1$-tree and\/ $S$ is a  subset of\/
 $\omega_1$ consisting of limit ordinals.  We say that\/ $T$ is
$S$-$*$-special iff there is a function $f$ mapping $T_S$ into $\omega_1$ such that
$(\forall x\in T_S)\allowbreak(f(x)<{\rm rk}(x))$ and whenever\/
$\{x,y,z\}\subseteq T_S$ and\/ $f(x)=f(y)=f(z)$ and\/ $x<y$ and\/ $x<z$ then\/
$y$ and\/ $z$ are comparable.

\proclaim Lemma 2.6.  Suppose\/ $T$ is an\/ $\omega_1$-tree and\/
$S_1$ and\/ $S_2$ are subsets of\/ $\omega_1$ consisting of limit ordinals and the 
symmetric difference\/
$S_1\Delta S_2=(S_1-S_2)\cup(S_2-S_1)$ is nonstationary.  
Then\/ $T$ is\/ $S_1$-$*$-special in the sense of Definition 2.5 iff\/
 $T$ is\/ $S_2$-$*$-special in the sense of Definition 2.5.  

Proof:  It suffices to
show that whenever $C$ is a closed unbounded set and $T$ is
$(S\cap C)$-$*$-special in the sense of Definition 2.5, then $T$ is $S$-$*$-special
in the sense of Definition 2.5.
Let $f$ mapping $T_{S\cap C}$ into $\omega_1$ be a specializing function.
For $x\in T_{S\cap C}$, let $\gamma_x$ be a limit ordinal (or zero) and $n_x$ an integer such that
$f(x)=\gamma_x+n_x$.  Because $S-C$ is a non-stationary set of limit ordinals,
 we may take $h$ to be a one-to-one function from
$S-C$ into $\omega_1$ such that $(\forall\alpha\in S-C)\allowbreak(h(\alpha)<\alpha)$.
We may assume that the range of $h$ consists only of odd ordinals.
Define $g$ such that for $x\in T_{S\cap C}$ we have that
$g(x)=\gamma_x+2n_x$, and for $x\in T_{S-C}$ we have $g(x)=h({\rm rk}(x))$.  Clearly $g$ demonstrates that $T$ is $S$-$*$-special in the sense
 of Definition 2.5.
The Lemma is established.

\proclaim Lemma 2.7.  Suppose\/ $T$ is an Aronszajn tree and\/ $E\subseteq\omega_1$
is uncountable.  Then\/ $T$ is\/ $E$-special in the sense of Definition 2.1 iff
there is\/ $g$ mapping\/ $T_E$ into\/ $\omega$ such that whenever\/ $x<y$ are in
$T_E$ then $g(x)\ne g(y)$. Furthermore, $T$ is $E$-special in the sense of
Definition 2.4 iff there is $h$ mapping $T_E$ into $\omega$ such that
whenever $x<y$ are in $T_E$ and\/ $z\in T_E$ and\/ $x<z$ and\/
$g(x)=g(y)=g(z)$ then $y$ is comparable with $z$.

Proof: The ``only if'' direction of the first assertion is evident by 
considering the composition of
a specializing function with a one-to-one mapping of ${\bf Q}$ into~$\omega$.  
For the ``if'' direction, suppose $g$ is given.  
Build $\langle f_n\,\colon n\in\omega\rangle$ by recursion such that 
$\dom(f_n)=\{x\in T_E\,\colon\allowbreak g(x)\leq n\}$ and 
${\rm range}(f_n)$ is a finite subset of ${\bf Q}$ and whenever $x<y$ are in ${\rm dom}
(f_n)$ then $f_n(x)<f_n(y)$.  
There is no difficulty in doing this. Clearly
$\bigcup\{f_n\,\colon\allowbreak n\in\omega\}$ is an $E$-specializing function in the sense of
Definition 2.1.  The first assertion is established.

The ``only if'' direction of the second assertion is again easy to see by considering
the composition of an $E$-specializing function with a mapping from
 ${\bf Q}$ into~$\omega$.  For the ``if'' direction, given $g$
  as in the statement of the assertion, then for
every $m\in\omega$ let ${\cal I}_m$ be the set of minimal elements of
$g^{-1}(m)$, and for every $x\in{\cal I}_m$ use the fact that
$T$ is Aronszajn to choose $\langle y_{m,x,i}\,\colon\allowbreak
i\in\omega\rangle$ an enumeration of $\{y\in T_E\,\colon\allowbreak x\leq y$ and $g(y)=m\}$. 
Build $\langle f_n\,\colon\allowbreak n\in\omega\rangle$ by recursion such that
 ${\rm dom}(f_n)=\{y\in T_E\,\colon \allowbreak
(\exists m\leq n)\allowbreak(\exists x\in{\cal I}_m)\allowbreak(\exists i\leq n)\allowbreak
(y=y_{m,x,i})\}$,
and ${\rm range}(f_n)$ is a finite subset of ${\bf Q}$ and
whenever $u<v$ and $u<w$ are all in ${\rm dom}(f_n)$ and $f_n(u)=f_n(v)=f_n(w)$ then 
$v$ is comparable with $w$.  There is again no difficulty in doing this.
Clearly $\bigcup\{f_n\,\colon\allowbreak
n\in\omega\}$ is an $E$-specializing function in the sense of Definition 2.2.
The Lemma is established.

\proclaim Lemma 2.8.  Suppose $T$ is Aronszajn and\/ $E\subseteq\omega_1$
 is unbounded.
 Then $T$ is $E$-special in the sense of Definition 2.1 iff $T$ is
$E$-special in the sense of Definition 2.4, and for $S\subseteq\omega_1$ consisting 
of limit ordinals we have
that\/ $T$ is
$S$-$*$-special in the sense of Definition 2.2 iff\/ $T$ is\/ 
$S$-$*$-special in the sense of Definition 2.5.

Proof: It is clear that if $T$ is $E$-special in the sense of Definition 2.1
then $T$ is $E$-special in the sense of Definition 2.4.  Suppose, therefore, that $T$ is
$E$-special in the sense of Definition 2.4.  Fix $f$ mapping
$T_E$ into ${\bf Q}$ as in Definition 2.4.
  Let $p$ be a one-to-one mapping from
${\bf Q}\,{\sf X}\,\omega$ into $\omega$.  For each $r\in{\bf Q}$ let 
${\cal I}_r$ be the set of all minimal elements of $f^{-1}(r)$. 
Using the fact that $T$ is Aronszajn, for each $x\in {\cal I}_r$ we may let
$\langle t_{r,x,k}\,\colon k\in\omega\rangle$ enumerate
 $\{y\in T_E\,\colon x\leq y$ and $f(y)=r\}$.  For every $z\in T_E$ let 
$h(z)=p(r,k)$ for the unique $r$ and $x$ and $k$ such that $z=t_{r,x,k}$.
It is clear that $h$ maps $T_E$ into $\omega$ and whenever $x<y$ are in
$T_E$ then $h(x)\ne h(y)$.  By Lemma 2.7 we have that $T$ is $E$-special
 in the
sense of Definition 2.1.

Now suppose that $S\subseteq\omega_1$ consists of limit ordinals.
Clearly if $T$ is $S$-$*$-special in the sense of Definition 2.2
 then $T$ is $S$-$*$-special in the sense of Definition 2.5.
So, suppose that $f$ is a function that $S$-$*$-specializes $T$
 in the sense of Definition 2.5.
By Lemma 2.6, we may assume that for every $\alpha\in S$ we have that 
 $\alpha=\alpha^\omega$ (ordinal arithmetic).
For each $\gamma\in\omega_1$ let ${\cal I}_\gamma$
 be the set of minimal elements
of $f^{-1}(\gamma)$.  For each $\gamma<\omega_1$ and
 $x\in {\cal I}_\gamma$ let
 $\langle t_{\gamma,x,m}\,\colon\allowbreak m\in\omega\rangle$
enumerate $\{y\in T_S\,\colon\allowbreak x\leq y$ and $f(y)=\gamma\}$.
 For $\gamma\in\omega_1$ and $x\in {\cal I}_\gamma$ and $m\in\omega$,
 set
$$g(t_{\gamma,x,m})=\omega\cdot\gamma+m$$
Using the fact that ${\rm rk}(t_{\gamma,x,m})\geq
{\rm rk}(x)>\gamma$ and ${\rm rk}(t_{\gamma,x,m})\in S$ we have that
${\rm rk}(t_{\gamma,i,m})\geq\gamma^\omega>g(t_{\gamma,x,m})$.
It is straightforward to check that $g$ is a function
 that $S$-$*$-specializes $T$ in the
sense of Definition 2.2.    

The Lemma is established.

\proclaim Lemma 2.9.  Suppose\/ $T$ is an\/ $\omega_1$-tree and either\/ $T$ is\/
 $E$-special for some 
unbounded\/ $E\subseteq\omega_1$ or\/ $T$ is\/ $S$-$*$-special for 
some stationary\/\ $S$ consisting of
limit ordinals.  Then\/ $T$ is neither Suslin nor Kurepa.

Proof:  By Lemmas 2.3 and 2.8 we have that $T$ is not Suslin.

Suppose that $T$ is $E$-special.  Let $f$ be a specializing function 
mapping $T_E$ into ${\bf Q}$.  
  For every uncountable branch $b$, 
choose $r_b\in{\bf Q}$ such that $\{y\in b\cap T_E\,\colon\allowbreak f(y)=r_b\}$ 
is uncountable, and let $t_b={\rm min}\{y\in b\cap T_E\,\colon\allowbreak f(y)=r_b\}$.
The function taking $b$ to $t_b$ is a one-to-one mapping from the set of uncountable branches into $T$.  Hence the
number of uncountable branches is at most $\aleph_1$.

Now assume that $S$ is a stationary set of countable limit ordinals 
and $T$ is $S$-$*$-special, and assume that $f$ is a function which $S$-$*$-specializes $T$.
    For every uncountable branch $b
\subseteq T$, use Fodor's Theorem to choose $\gamma_b\in\omega_1$ such that $\{y\in b\cap T_S\,\colon\allowbreak f(y)=\gamma_b\}$ is uncountable, and let $t_b={\rm min}\{y\in b\cap T_S\,\colon\allowbreak f(y)=\gamma_b\}$.
The function taking $b$ to $t_b$ is a one-to-one mapping from the set of uncountable branches into $T$.  Hence the
number of uncountable branches is at most $\aleph_1$.  The Lemma is established.

\section{$(T,S)$-$\#$-preserving forcings}

In this section we introduce the preservation property that will be
used in the main constructions, and we establish some technical properties.

\proclaim Definition 3.1.  Suppose 
$T$ is an $\omega_1$-tree and\/ $\lambda$ is a sufficiently large regular cardinal and
$N$ is a countable elementary substructure of $H_\lambda$ and
$T\in N$ and $x\in T$.  We say that\/ $x$ is $(T,N)$-$\#$-generic iff 
${\rm rk}(x)=\omega_1\cap N$ and for every $A\in N$
such that $A\subseteq T$
 we have\/ {\rm $(\exists y< x)\allowbreak(y\in A$ or
$(\forall z\geq y)\allowbreak(z\notin A))$.}

\proclaim Definition 3.2.  Suppose $P$ is a forcing and $T$ is an $\omega_1$-tree and
$S\subseteq\omega_1$ and $\lambda$ is a sufficiently large regular cardinal and 
$N$ is a countable elementary substructure of
$H_\lambda$ and $\{P,T, S\}\in N$ and $q\in P$.  We say that 
$q$ is $(N,P,S,T)$-$\#$-preserving iff
$q$ is $(N,P)$-generic and either $\omega_1\cap N\in S$ or
for every $x\in T$ such that 
$x$ is $(T,N)$-$\#$-generic
 and every  $P$-name 
$A$ from $N$ such that\/ {\rm ${\bf 1}\forces``A\subseteq
T$,''} we have that\/ {\rm
$q\forces``(\exists y< x)\allowbreak(
y\in A$ or
 $(\forall z\geq y)\allowbreak(z\notin A))$.''}

\proclaim Lemma 3.3. Whenever $p$ is $(N,P,S,T)$-$\#$-preserving and $q\leq p$, then
$q$ is $(N,P,S,T)$-$\#$-preserving.  Also, whenever $p$ is $(N,P,S,T)$-$\#$-preserving and
$S\subseteq S'$ then $p$ is $(N,P,S',T)$-$\#$-preserving. 

Proof: Obvious.

\proclaim Lemma 3.4.  Suppose $q$ is $(N,P,S,T)$-$\#$-preserving and 
$\omega_1\cap N\notin S$ and $x$ is
$(T,N)$-$\#$-generic. Then
{\rm $q\forces``x$ is
$(T,N[G_P])$-$\#$-generic.''}

Proof: Necessarily $q$ is $N$-generic, so $q\forces``{\rm rk}(x)=\omega_1\cap N=\omega_1\cap
N[G_P]$ and $T$ is an $\omega_1$-tree.'' Now suppose that
$q'\leq q$ and 
$q'\forces``
A\in N[G_P]$ and
$A\subseteq T$ and $(\forall y<x)\allowbreak(y\notin A)$.''
Because $q'\forces``A\in N[G_P]$'' we may take
$r\leq q'$ and $A'$ a $P$-name in $N$ such that
$r\forces``A'=A$.''
We may replace $A'$ by the $P$-name $A^*$ in $N$ characterized
by ${\bf 1}\forces``A^*=A'$ if $A'\subseteq
T$ and $A^*=\emptyset$ otherwise.''
Because $r$ is $(N,P,S,T)$-$\#$-preserving, we have
$r\forces``(\exists y<x)\allowbreak(\forall z\geq y)
\allowbreak(z\notin A^*)$.''
The Lemma is established.

\proclaim Lemma 3.5.  Suppose $p$ is $(N,P,S,T)$-$\#$-preserving and\/ {\rm
$p\forces``\dot q$  is
$(N[G_P],\dot Q,\allowbreak S,T)$-$\#$-preserving.''} Then
$(p,\dot q)$ is $(N,P*\dot Q,S,T)$-$\#$-preserving.

Proof:  If $\omega_1\cap N\in S$, then the Lemma
follows from the well-known fact 
that if $p$ is $N$-generic and $p\forces``\dot q$ is $N[G_P]$-generic,'' then
$(p,\dot q)$ is $N$-generic.
So suppose that $\omega_1\cap N\notin S$  
 and $A\in N$
is a $P*\dot Q$-name for a subset of $T$  and $x$ is $(T,N)$-$\#$-generic and
$(p_1,\dot q_1)\leq(p,\dot q)$ and
$(p_1, \dot q_1)\forces``(\forall y<x)\allowbreak
(y\notin A)$.''
Fix $ \tilde A\in N$ a $P$-name such that
${\bf 1}\forces_P`` \tilde A$ is a $\dot Q$-name and ${\bf 1}\forces_{\dot Q}` \tilde A=
 A$.'\thinspace''  Because $p$ is $N$-generic, we have that
$p\forces``T$ is an $\omega_1$-tree.''

By Lemmas 3.3 and 3.4, we have $p_1\forces``x$ is $(T,N[G_P])$-$\#$-generic.'' 
Because $p_1$ is $N$-generic 
we also have that
$p_1\forces``\omega_1\cap N[G_P]=\omega_1\cap N\notin S$.'' 
Hence using the fact that $p_1\forces``\dot q_1$ is
 $(N[G_P],\dot Q,S,T)$-$\#$-preserving
and $\dot q_1\forces`(\forall y<x)\allowbreak
(y\notin \tilde A)$,'\thinspace'' 
we have that
$p_1\forces``\dot q_1\forces`(\exists y<x)\allowbreak(\forall z\geq y)
\allowbreak(z\notin \tilde A)$.'\thinspace''
So there is $(p_2,\dot q_2)\leq(p_1,\dot q_1)$ and $y<x$ such that
$p_2\forces``\dot q_2\forces`(\forall z\geq y)\allowbreak(z\notin \tilde A)$.'\thinspace''
We have $(p_2,\dot q_2)\forces``(\forall z\geq y)\allowbreak
(z\notin A)$.'' The Lemma is established.

The following is the key Definition of this paper, in which we
isolate the preservation property that we use to maintain
non-$E$-specialness of an appropriately chosen
 Suslin tree  of the ground model.  This
Definition is analogous to [PIF, Definition IX.4.5], [JSL, Definition 5].

\proclaim Definition 3.6.  Suppose $T$ is  an $\omega_1$-tree and $S\subseteq\omega_1$  and 
$P$ is a poset.  We say that $P$ is $(T,S)$-$\#$-preserving iff
whenever $\lambda$ is a sufficiently large regular cardinal and $N$ is a countable 
elementary substructure of $H_\lambda$ and
$\{T,S,P\}\in N$ and $p\in P\cap N$ then there is 
$q\leq p$ such that $q$ is $(N,P,S,T)$-$\#$-preserving.

In the following Definition, we specify three different 
ways of collapsing a stationary co-stationary subset of $\omega_1$.  
These are well-known (although the third poset
 is less well-known than it deserves to be).

\proclaim Definition 3.7. Suppose $S\subseteq\omega_1$.  The poset $CU(S)$ 
is the set of closed, bounded subsets of $S$ ordered by reverse end-extension.  
The poset $CU^*(S)$ consists of pairs $\langle 
\sigma,C\rangle$ such that $\sigma $ is a countable closed subset of $S$ and $C$ 
is a closed unbounded subset of
$\omega_1$, ordered by $\langle\sigma_1,C_1\rangle\leq\langle\sigma_2,C_2\rangle$ 
iff $\sigma_1$ end-extends
$\sigma_2$ and $C_1\subseteq C_2$ and $\sigma_1\subseteq\sigma_2\cup C_2$.  
The poset $CU^{**}(S)$ consists of
all finite sets $F$ of intervals $[\alpha,\beta]$ such that the elements of $F$ 
are disjoint, and for every
$[\alpha,\beta]\in F$ we have that $\alpha$ is either a successor ordinal or zero 
or an element of $S$, ordered by $F_1\leq F_2$ iff $F_1\supseteq F_2$.

\proclaim Lemma 3.8. Suppose $S\subseteq\omega_1$.  Suppose $P$ is one of $CU(S)$ 
or $CU^*(S)$ or $CU^{**}(S)$. Then in $V[G_P]$ we have that $\omega_1-S$ is 
non-stationary, and if $S$ is stationary then $\omega_1$ is preserved
(in fact,
$P$ is $S$-proper).

Proof: The only possibly unclear case is handled by the observation that 
if $P=CU^{**}(S)$ then in $V[G_P]$, we have that
$\{\alpha\,\colon\allowbreak(\exists F\in G_P)\allowbreak(\exists\beta<\omega_1)
\allowbreak([\alpha,\beta]\in F$ and $\alpha$ is a limit ordinal$)\}$ is 
a closed unbounded subset of $S$.  The $S$-properness of
$CU^{**}(S)$ is demonstrated in the proof of Lemma 3.9.

The ``case 2'' part of the proof of Lemma 3.9
 recalls the proof of [Sh, Lemma IX.4.6]. Lemma 3.9 is analogous to
[PIF, XXX] and [JSL, Lemma 20].

\proclaim Lemma 3.9.  Suppose $P=CU(S)$ or $P=CU^*(S)$ or $P=CU^{**}(S)$ and
$\lambda$ is a sufficiently large regular cardinal and
$N$ is a countable elementary substructure of $H_\lambda$ and $\{T,S,P\}\in N$ and
$\omega_1\cap N\in S$  and $p\in P\cap N$.  Then there is $q\leq p$
such that $q$ is $(N,P,\omega_1-S,T)$-$\#$-preserving.

Proof: Let $\delta=\omega_1\cap N$.  

Case 1: $P=CU(S)$ or $P=CU^*(S)$.

Let $\langle (x_n,A_n)\,\colon n\in\omega\rangle$ enumerate the set of all
$\langle x,A\rangle$ such that $A\in N$ is a $P$-name for a subset of $T$ and $x$ is
$(T,N)$-$\#$-generic (if there are no such $x$ then ignore
requirement (3) below).  
  Let $\langle D_n\,\colon\allowbreak
n\in\omega\rangle$ list the set of all $D\in N$ such that $D\subseteq P$ is open dense.

Claim: There is a sequence $\langle p_n\,\colon n\in\omega\rangle$ such that $p_0=p$ and
for all $n\in\omega$ we have that each of the following holds:

(1) $p_{n+1}\leq p_n$

(2) $p_{n+1}\in D_n\cap N$

(3) either $p_{n}\forces``(\exists y< x_n)(\forall z\geq y)(z\notin A_n)$'' or for some $y<x_n$ we have $p_{n+1}\forces``
y\in A_n$.''

Proof of Claim:  
Given $p_n$, take $p'_n\leq p_n$ such that
$p'_n\in D_n\cap N$.  Let $Y=\{y\in T\,\colon\allowbreak
p'_n\notforces``y\notin A_n$''$\}$.

If $(\exists y<x_n)\allowbreak(y\in Y)$ then we may take $p_{n+1}\leq p'_n$
 and $y<x_n$ such that
$p_{n+1}\forces``y\in A_n$.''  We may assume $p_{n+1}\in N$, and hence the
second disjunct in requirement (3) holds.

If instead $(\forall y<x_n)\allowbreak(y\notin Y)$, then because
 $x_n$ is $(T,N)$-$\#$-generic we have that there is some $y<x_n$
  such that $(\forall z\geq y)\allowbreak
(z\notin Y)$.  Hence the first disjunct of requirement (3) holds.

If $P=CU(S)$ then let $q=\bigcup\{p_n\,\colon\allowbreak n\in\omega\}\cup\{\delta\}$, whereas
if $P=CU^*(S)$ then let
$q=\langle\bigcup\{\sigma_n\,\colon\allowbreak n\in\omega\}\cup\{\delta\},\allowbreak
\bigcap\{C_n\,\colon\allowbreak n\in\omega\}\rangle$ where $p_n=\langle\sigma_n,C_n\rangle$ for every
$n\in\omega$.  Because $\delta\in S$ we have
$q\in P$.  Clearly $q$ is as required.

Case 2:  $P=CU^{**}(S)$.

Let $\delta^*={\rm sup}\{f(\delta)+1\,\colon\allowbreak f\in N$ is a function
and $f(\delta)\in\omega_1\}$.  Let $q=p\cup\{[\delta',\delta']\}$, where $\delta'\geq\delta^*$ is not
a limit ordinal outside of $S$ (hence $q\in P$).
We show that $q$ is $(N,P,\omega_1-S,T)$-$\#$-preserving.

First we show that $q$ is $(N,P)$-generic. Given $D\in N$ a dense open
 subset
of $P$, and given $q^*\leq q$,
we find $r\leq q^*$ such that $r$ is below some element of $D\cap N$.
Choose $r'\leq q^*$ such that $r'\in D$.  We  have 

\medskip

\centerline{$N\models``(\exists p^*\leq (r'\cap N))(p^*\in D)$''}

\medskip

Choose $p^*\in  N$ to be a witness. Set $r=r'\cup p^*$.  Clearly
$r\in P$ and $r$ is as required.

Now suppose, towards a contradiction, that $x$ is $(T,N)$-$\#$-generic and
 $A\in N$ is a $P$-name for a subset of
$T$, and $q'\leq q$ and $q'\forces``
(\forall y< x)\allowbreak(y\notin A$
and $(\exists z\geq y)\allowbreak
(z\in A))$.''

Let $\alpha={\rm sup}(\bigcup(q'\cap N))$.  In other words, the
``largest'' interval in
$q'\cap N$ is $[\gamma,\alpha]$ for some $\gamma$.

For $p_1$ and $p_2$ in $P$, define $p_1\leq^* p_2$ iff
there is some $\beta$ such that 
$p_2=\{[\eta,\gamma]\in p_1\,\colon\allowbreak\gamma\leq\beta\}$.
Essentially, $p_1\leq^* p_2$ iff $p_1$ ``end-extends'' $p_2$.

Let $R=\{y\in T\,\colon\allowbreak
{\rm rk}(y)>\alpha\}$ and for all $y\in R$ let
 $J(y)=\{\gamma<\omega_1\,\colon\allowbreak
(\exists\alpha^*<\omega_1)\allowbreak(\gamma\leq\alpha^*$ and
$\gamma$ is not a limit ordinal outside of $S$
and $(\exists\tilde q\leq^* (q'\cap N
\cup\{[\gamma,\alpha^*]\}))\allowbreak
(\tilde q\forces``y\notin A$''$))\}$.
Let $F$ be the function with domain equal to $\{y\in R\,\colon\allowbreak
J(y)\ne\emptyset\}$ characterized by $(\forall y\in {\rm dom}(F))
\allowbreak(F(y)={\rm sup}(J(y)))$.
Let $A^*=\{y\in R\,\colon J(y)\ne\emptyset$ and $ F(y)=\omega_1\}$.

Because $x$ is $(T,N)$-$\#$-generic, we may fix $y<x$ such that
either $y\notin A^*$ or $(\forall z\geq y)\allowbreak
(z\in A^*)$.

Case 1.  $y\notin A^*$

Claim: $(q'\cap N)\notforces``y\notin A$.''

Suppose instead that $(q'\cap N)\forces``y\notin A$.'' We have
$q'\cap N\cup\{[\gamma+1,\gamma+1]\}$ witnesses $\gamma+1\in J(y)$ for every
countable $\gamma\geq\alpha$, hence $F(y)=\omega_1$, contradicting the fact that 
$y\notin A^*$. The Claim is established.

By the Claim  we may take $q^+\leq (q'\cap N)$ such that $q^+\forces``y\in A$'' and
$q^+\in N$.  Clearly we have that $(q^+\cup q')\in P$. But $q^+
\forces``y\in A$'' and $q'\forces``y\notin A$.'' This is impossible.

Case 2.  $(\forall z\geq y)(z\in A^*)$

We have $q'\forces``(\exists z\geq y)(z\in A)$.'' Choose $q^+\leq q'$
and $z\geq y$ such that $q^+\forces``z\in A$.''
Fix $\gamma$ a countable ordinal greater than ${\rm sup}(\bigcup q^+)$.
Because $z\in A^*$ we know that 
$J(z)$ is not empty. Furthermore, $F(z)=\omega_1$, so we may take
$\gamma^*\in J(z)$ such that $\gamma\leq\gamma^*$. We may $\alpha^*\geq \gamma^*$
and $\tilde q\leq^*
(q'\cap N\cup\{[\gamma^*,\alpha^*]\}$ such that $\tilde q\forces``z\notin
A$.'' 
Clearly $(q^+\cup\tilde q)\in P$.  We have
$(q^+\cup\tilde q)\forces``z\in A$ and $z\notin A$.'' This is impossible, hence
the Lemma is established.

The following Lemma is analogous to [JSL, Lemma 21].

\proclaim Lemma 3.10.  Suppose $P$ is a poset and $T$ is an
$\omega_1$-tree and $S\subseteq\omega_1$.
Suppose $\lambda$ is a sufficiently large regular cardinal and
 $N$ is a countable elementary substructure of $H_\lambda$ containing
 $\{P,T,S\}$.  Suppose
$p\in P$ is $(N,P,S,T)$-$\#$-preserving and\/ $A$ is a $P$-name in $N$ for a $Q$-name in
$N[G_P]$ that names a
subset of $T$ and\/ {\rm
 $p\forces``\dot Q$ is $(T,S)$-$\#$-preserving and $\dot Q\in N[G_P]$ and
$q\in\dot Q\cap N[G_P]$''}
  and $x\in T$ and
$x$ is $(T,N)$-$\#$-generic and $\omega_1\cap N\notin S$.  Then there is
 a $P$-name $r$ such
that\/ {\rm $p\forces``r\leq q$ and $r\in N[G_P]$''} and\/
{\rm $(p,r)\forces``(\exists y<x)\allowbreak(y\in A$ or
$(\forall z\geq y)\allowbreak(z\notin A))$.''}

Proof: Let $D=\{p'\leq p\,\colon\allowbreak 
p'\forces``(\exists y<x)\allowbreak(
q\forces`(\forall z\geq y)\allowbreak
(z\notin A)$')'' or 
$p'\forces``(\exists y<x)\allowbreak(\exists q'\leq q)
\allowbreak(q'\in N[G_P]$ and
$q'\forces`y\in A$')''$\}$.

Claim 1.  $D$ is dense below $p$.

Proof: Suppose $p^+\leq p$. Because $p\forces``q\in N[G_P]$,'' we may take
$\tilde p\leq p^+$ and $q^*$ a $P$-name in $N$ such that
$\tilde p\forces``q^*=q$.'' Take $B$ to be a $P$-name in $N$ characterized by
${\bf 1}\forces``B=\{y\in T\,\colon\allowbreak q^*\notforces`y\notin A$'$\}$.'' 

By Lemma 3.4 we have
 $\tilde p\forces``x$ is $(T,N[G_P])$-$\#$-generic,'' and therefore
  we can take
$p_1\leq \tilde p$ and $y<x$ such that either
$p_1\forces``y\in B$'' or
$p_1\forces``(\forall z\geq y)\allowbreak(z\notin B)$.''

If we have $p_1\forces``y\in B$,'' then  $ p_1$ witnesses the second
disjunct in the
definition of $D$ and we are done. If instead 
$p_1\forces``(\forall z\geq y)\allowbreak(z\notin B)$,''
we have $p_1\forces``(\forall z\geq y)(q^*\forces`z\notin A$'),''
and thus $p_1$ witnesses the first disjunct of the definition of $D$.
In either case, the Claim is established.

We now define a function $f$ with domain $D$ as follows.  If $p'\in D$ and $p'\forces``
(\exists y<x)(
q\forces`(\forall z\geq y)\allowbreak(z\notin A)$'),''
 then we let $f(p')=q$. 
If instead $p'\forces``(\exists q'\leq q)\allowbreak(\exists y<x)\allowbreak
(q'\in N[G_P]$ and $q'\forces`y\in A$'))'' then we choose
some such $q'$, and 
set $f(p')=q'$.
  Let ${\cal J}$ be a maximal 
antichain of $P$ such
 that ${\cal J}\subseteq D$.
Let $r$ be a $P$-name such that for every $p'\in{\cal J}$ we have
$p'\forces``r=f(p')$.''

By  Claim 1 we clearly have that $p\forces``r\in N[G_P]$ and $r\leq q$.''

Claim 2. $p\forces``r\forces`(\exists y<x)\allowbreak
(y\in A$ or $(\forall z\geq y)\allowbreak
(z\notin A))$.'\thinspace''

Proof: Suppose $p_1\leq p$. Take $p'\in{\cal J}$ and $p_2\leq p_1$ such that
$p_2\leq p'$.

Case 1: $p'\forces``(\exists y<x)\allowbreak(q\forces`(\forall z\geq y)\allowbreak(z\notin
A)$').''

Clearly $p'\forces``r\forces`(\exists y<x)\allowbreak(\forall z\geq y)
\allowbreak(z\notin A)$.'\thinspace''

Case 2: Otherwise.

Because Case 1 fails and $p'\in D$ we have by choice of $r$ 
that $p'\forces``(\exists y<x)\allowbreak(r\forces`y\in A$').''

In either case, we have $p_2\forces``r\forces`(\exists y<x)\allowbreak
(y\in A$ or $(\forall z\geq y)(z\notin A))$'\thinspace''

The Claim is established, and the Lemma is proved.

\proclaim Lemma 3.11. Suppose $x$ is $(T,N)$-$\#$-generic and $A\in N$
and $A$ is an antichain of $T$ and $T\in N$. Then $x\notin A$.

Proof: Suppose $x\in A$. Then $(\forall y<x)\allowbreak
(y\notin A)$. Hence $(\exists y<x)\allowbreak
(\forall z\geq y)\allowbreak(z\notin A)$. Hence $x\notin A$.
The Lemma is established.

\section{Sharply Suslin trees}

\proclaim Definition 4.1.  Suppose $T$ is a
Suslin tree.  We say that $T$ is
sharply Suslin iff 
 for every sufficiently large regular cardinal $\lambda$  and every 
countable elementary substructure $M$ of $H_{\lambda^+}$ we have that
 there is
a closed unbounded $C\subseteq\omega_1$ such that for every
$N\in M$ such that $N$
is a countable elementary substructure of $H_\lambda$ containing
$T$ and $\omega_1\cap N\in C$, and
 every $x\in T_{\omega_1\cap N}$, we have that $x$
  is $(T,N)$-$\#$-generic.

\proclaim Definition 4.2. $\diamondsuit^*$ is the following principle:
there is a sequence
$\langle S_\alpha\,\colon\allowbreak\alpha<\omega_1\rangle$ such that
$S_\alpha$ is a countable subset of ${\cal P}(\alpha)$ and for every
 $X\subseteq\omega_1$
 there is
a closed unbounded $C\subseteq\omega_1$ such that for
 every $\alpha\in C$ we have
$X\cap\alpha\in S_\alpha$.

\proclaim Lemma 4.3. Suppose $V=L$. Then $\diamondsuit^*$ holds.

Proof: See [Devlin, Theorem III.3.5].

\proclaim Lemma 4.4.  Suppose $\diamondsuit^*$ holds.
Then there is a sharply Suslin tree.

Proof: Let 
$\langle S_\alpha\,\colon\allowbreak \alpha<\omega_1\rangle$ be a
$\diamondsuit^*$-sequence.  Because $\diamondsuit^*$ implies $\diamondsuit$,
we may also fix a $\diamondsuit$-sequence
$\langle Z_\alpha\,\colon\allowbreak\alpha<\omega_1
\rangle$.

Given $\lambda$ and $M$ as in Definition 4.1,
let $X\subseteq\omega_1$ code $\{x\cap\omega_1\,\colon\allowbreak
x\in M\}$.  For example, we may let
$\langle\theta_i\,\colon\allowbreak i\in\omega\rangle$ list
$M$, and let $X=\{\omega\alpha+i\,\colon\allowbreak\alpha\in\theta_i\}$.
Let $C$ be a closed unbounded subset of
$\{\alpha<\omega_1\,\colon\allowbreak
X\cap\alpha\in S_\alpha$ and $\alpha$ is indecomposable$\}$.
Build $T$ recursively such that whenever $\beta<\omega_1$ is
an indecomposable ordinal then we have
 $T_{<\beta}=\beta$ and we build $T_\beta$ as follows.

Let $\langle\beta_n\,\colon\allowbreak n\in\omega\rangle$ be an increasing
sequence of ordinals cofinal in $\beta$.  Let $\langle B^\beta_i\,\colon
\allowbreak i\in\omega\rangle$ list $S_\beta$. Let
$\langle A^\beta_i\,\colon\allowbreak i\in\omega\rangle$ list
$\{\{\alpha<\beta\,\colon\allowbreak
\omega\alpha+i\in B^\beta_k\}\,\colon\allowbreak i\in\omega$ and
$k\in\omega\}$.  Thus for each $x\in M$, if $\beta\in C$ then we have that
$x\cap\beta$ is equal to $A^\beta_i$ for some $i\in\omega$.

Build $\langle (A')^\beta_i\,\colon\allowbreak i\in\omega\rangle$ such that
 for every $n\in\omega$ we have
 each of the following:

(1)  $(A')^\beta_n$ is an antichain of $T_{<\beta}$,

(2) $(A')^\beta_n$ is  predense above $A^\beta_n$, i.e.,
$(\forall x\in A^\beta_n)\allowbreak
(\forall y\geq x)\allowbreak
(\exists z\in(A')^\beta_n)\allowbreak
(z$ is comparable with $y)$,

(3) for every $i\leq n$ and every $y'\in (A')^\beta_i$ and every
$y\in (A')^\beta_n$ we have either  $y'\leq y$ or
$y'$ is incomparable with  $y$,

(4) for every $y\in(A')^\beta_n$ we have ${\rm rk}(y)\geq\beta_n$.

There is no problem in doing this.

Now construct $T_\beta$ such that

(1) for every $x\in T_\beta$ and every $n\in\omega$ there is
$y<x$ such that either $y\in (A')^\beta_n$ or $(\forall z\geq y)\allowbreak
(z\notin (A')^\beta_n)$, and

(2) for every $y\in T_{<\beta}$ there is $x\in T_\beta$ such that
$y<x$,

(3) for every $x\in T_\beta$ there is $y<x$ such that either
$y\in Z_\beta$ or $(\forall z\geq y)\allowbreak
(z\notin Z_\beta)$.

There is no problem in this.

It is easy to see that for every $x\in T_\beta$ and every $n\in\omega$ there
is $y<x$ such that either $y\in A^\beta_n$ or $(\forall z\geq y)\allowbreak
(z\notin A^\beta_n)$.  Also it is easy to see that the tree $T$ that is
constructed in this way is a Suslin tree.

It is easy to see that $C$ is the required witness to the assertion that
$\lambda$ and $M$ do not constitute a counterexample to the fact that
$T$ is sharply Suslin.

The Lemma is established.

\proclaim Lemma 4.5.  Suppose $T$ is a sharply Suslin
tree and $S\subseteq\omega_1$ is
co-stationary and $P$ is $(T,S)$-$\#$-preserving and\/ {\rm ${\bf 1}\forces``T$ is
Aronszajn.''} Then\/ {\rm
${\bf 1}\forces``(\forall E\subseteq\omega_1$ unbounded$)\allowbreak
(T$ is not $E$-special).''}

Proof: Suppose, towards a contradiction, that $E$ and $f$ are $P$-names
and $p\in P$ and $p\forces``E\subseteq\omega_1$ is unbounded and
$f$ is an $E$-specializing function for $T$ in the sense of Definition 2.1.''
  Take $\lambda$ a large enough regular cardinal and 
$M$ a countable
 elementary substructure of $H_{\lambda^+}$
containing $\{T,S,P,E,f,p\}$.  Take $C$ as in Definition 4.1 and fix
 $N\in M$ such that
 $N$ is a countable elementary substructure of $H_\lambda$ containing
$\{T,S,P,E,f,p\}$ and such that
 $\omega_1\cap N\in C$
and $\omega_1\cap N\notin S$.
 Fix $x\in T_{\omega_1\cap N}$.
Take $q\leq p$ such that $q$ is $(N,P,S,T)$-$\#$-preserving and, by a further strengthening
of $q$, we may take $r\in{\bf Q}$ such that $q\forces``f(z)=r$ for some $z\geq x$.''
Necessarily we have $q\forces``(\forall y<x)\allowbreak
(y\notin f^{-1}(r))$.''
Hence $q\forces``(\exists y<x)\allowbreak(\forall y'\geq y)\allowbreak(y'\notin f^{-1}(r))$.''
This contradicts the fact that $q\forces``z\in f^{-1}(r)$.''
 The Lemma is established.

\section{Iteration of $(T,S)$-$\#$-preserving forcings}

In this section we show that the property
 $``(T,S)$-$\#$-preserving'' is preserved
by countable support forcing iteration (and a bit more).  This is a
variant of [PIF XXXX], [JSL, Definition 22 and Lemmas 23 and 24].

\proclaim Definition 5.1.  Suppose
 $\langle P_\eta\,\colon\eta\leq\kappa\rangle$ is a countable support iteration of forcing. 
 We say that $P_\kappa$ is strictly $(T,S)$-$\#$-preserving iff
whenever $\lambda$ is a sufficiently large regular cardinal and $N$ 
is a countable elementary substructure of $H_\lambda$ and
$\{T,S,\langle P_\eta\,\colon\allowbreak\eta\leq\kappa\rangle\}\in N$ 
and $\alpha\in\kappa\cap N$ and $q\in P_\alpha$ is
$(N,P_\alpha,S,T)$-$\#$-preserving and\/ {\rm 
$q\forces``\dot p\in \dot P_{\alpha,\kappa}\cap N[G_{P_\alpha}]$,''}
 then there is
$r\in P_\kappa$ such that $r\restr\alpha=q$ and\/ {\rm 
$q\forces``r\restr[\alpha,\kappa)\leq\dot p$''} and
$r$ is $(N,P_{\kappa},S,T)$-$\#$-preserving
and $\supt(r)\subseteq \alpha\cup N$.

\proclaim Lemma 5.2.  Suppose $P_\kappa$ is strictly $(T,S)$-$\#$-preserving.  Then
$P_\kappa$ is $(T,S)$-$\#$-preserving.

Proof: Take $\alpha=0$ in Definition 5.1.

\proclaim Lemma 5.3.  The following are equivalent:

(1)  $P_\kappa$ is strictly $(T,S)$-$\#$-preserving,

(2) For some regular $\lambda>\omega_1$ such that $P_\kappa\in H_\lambda$ there is a closed unbounded
$C\subseteq[H_\lambda]^\omega$ such that whenever $N\in C$ and 
$\eta\in \kappa\cap N$ and
$q$ is $(N,P_\eta,S,T)$-$\#$-preserving and $q\forces``\dot p\in\dot P_{\eta,\kappa}\cap N[G_{P_\eta}]$''
then there is $r\in P_\kappa$ such that $r\restr\eta=q$ and 
$q\forces``r\restr[\eta,\kappa)\leq\dot p$''
and $r$ is $(N,P_{\kappa},S,T)$-$\#$-preserving and $\supt(r)\subseteq\eta\cup N$,

(3) For some regular $\lambda>2^{\aleph_1}$ such that the power set of $P_\kappa$ is an element of $H_\lambda$
we have that  whenever $N$ is a countable elementary substructure of $H_\lambda$
and $\{P_\kappa,T,S\}\in N$ and $\eta\in \kappa\cap N$ and
$q$ is $(N,P_\eta,S,T)$-$\#$-preserving and $q\forces``\dot p\in\dot P_{\eta,\kappa}\cap N[G_{P_\eta}]$''
then there is $r\in P_\kappa$ such that $r\restr\eta=q$ and 
$q\forces``r\restr[\eta,\kappa)\leq\dot p$''
and $r$ is $(N,P_{,\kappa},S,T)$-$\#$-preserving and $\supt(r)\subseteq	\eta\cup N$,

(4) For every regular $\lambda>\omega_1$ such that $P_\kappa\in H_\lambda$ there is a closed unbounded
$C\subseteq[H_\lambda]^\omega$ such that whenever $N\in C$ and $\eta\in \kappa\cap N$ and
$q$ is $(N,P_\eta,S,T)$-$\#$-preserving and $q\forces``\dot p\in\dot P_{\eta,\kappa}\cap N[G_{P_\eta}]$''
then there is $r\in P_\kappa$ such that $r\restr\eta=q$ and 
$q\forces``r\restr[\eta,\kappa)\leq\dot p$''
and $r$ is $(N,P_{,\kappa},S,T)$-$\#$-preserving and $\supt(r)\subseteq\eta\cup N$.

\medskip

Proof: (1) trivially implies (2) and (3), and (4) trivially implies (2). We show (2) implies (1). Fix $\lambda$ to be the least witness to (2), and suppose 
$\mu>2^{\aleph_1}$ is a regular cardinal such that the power set of $P_\kappa$ is in $H_\mu$.
Let $C_0=\{M\in[H_\lambda]^\omega\,\colon \{P_\kappa,T,S\}\in M$ and 
whenever $\eta\in\kappa\cap M$ and $q$ is
$(M,P_\eta,S,T)$-$\#$-preserving and $q\forces``\dot p\in\dot P_{\eta,\kappa}\cap M[G_{P_\eta}]$'' then there is
$r\in P_\kappa$ such that $r\restr\eta=q$ and 
$q\forces``r\restr[\eta,\kappa)\leq\dot p$'' and
$r$ is $(M,P_{\kappa},S,T)$-$\#$-preserving and $\supt(r)\subseteq\eta\cup M\}$.  Suppose $N$ is a countable elementary substructure of $H_\mu$
and $\{P_\kappa,T,S\}\in N$.  Then $C_0\in N$ and $\lambda\in N$, because
 $C_0$ and $\lambda$ are $\Delta_1$-definable from the parameters  $P_\kappa$, ${\cal P}(P_\kappa)$, ${\cal P}(\omega_1)$, $T$, and $S$.
Take $C\subseteq C_0$ such that $C$ is a closed unbounded subset of $[H_\lambda]^\omega$ and $C\in N$.
Let $\langle\theta_n\,\colon\allowbreak n\in\omega\rangle$ enumerate $N\cap H_\lambda$. By recursion, build
$\langle M_n\,\colon\allowbreak n\in\omega\rangle$ such that for every $n\in\omega$ we have
$M_n$ is a countable elementary substructure of $N\cap H_\lambda$ and $\{M_n,\theta_n\}\in M_{n+1}$ and
$M_n\in C$.  We have therefore that $N\cap H_\lambda=\bigcup\{M_n\,\colon\allowbreak n\in\omega\}\in C$.
Therefore, whenever $\eta\in\kappa\cap N$
 and $q$ is $(N,P_\eta,S,T)$-$\#$-preserving and 
$q\forces``\dot p\in \dot P_{\eta,\kappa}\cap N[G_{P_\eta}]$,'' then clearly $q$ is $(N\cap H_\lambda,P_\eta,S,T)$-$\#$-preserving
and $q\forces``\dot p\in (N\cap H_\lambda)[G_{P_\eta}]$,'' and therefore there 
is $r\in P_\kappa$ such that
$r\restr\eta=q$ and $q\forces``r\restr[\eta,\kappa)\leq\dot p$'' and
$r$ is $(N\cap H_\lambda,P_{\kappa},S,T)$-$\#$-preserving and $\supt(r)\subseteq
\eta\cup N$.  Clearly $r$ is $(N,P_{\kappa},S,T)$-$\#$-preserving.  This verifies that (1) holds.

We now show that (3) implies (4).  Given $\lambda$ as in (4), let $C=\{M\in [H_\lambda]^\omega\,\colon\allowbreak M$ is a countable elementary substructure of $H_\lambda$ and $\{P_\kappa,T,S\}\in M$ and there is some regular $\mu>2^{\aleph_1}$ and $N$ a countable elementary substructure of $H_\mu$ such that
$2^{\omega_1}\in N$ and ${\cal P}(P_\kappa)\in N$
and $M=N\cap H_\lambda\}$. Then $C$ witnesses that (4) holds.
 The Lemma is established.

\proclaim Theorem 5.4.  Suppose $T$ is sharply Suslin and $S\subseteq\omega_1$ and
$\langle P_\eta\,\colon\allowbreak\eta\leq\kappa\rangle$ is a countable 
support forcing iteration based on $\langle\dot Q_\eta\,\colon\allowbreak\eta<\kappa\rangle$.
Suppose for every $\eta<\kappa$ we have that ${\cal I}_\eta$ 
and $S_\eta$ are $P_\eta$-names.
  Suppose for every $\eta<\kappa$ we have either

(1) $\bfone\forces_{P_\eta}``\dot Q_\eta$ is $(T,S)$-$\#$-preserving,'' {\sl or}

(2) $\bfone\forces_{P_\eta}``{\cal I_\eta}$ is an antichain of $T$ and $S_\eta=\{{\rm rk}(x)\,\colon\allowbreak
x\in{{\cal I}_\eta}\}$ and $\dot Q_\eta$ is one of $CU(S\cup(\omega_1-S_\eta))$ or $CU^*(S\cup(\omega_1-S_\eta))$ or $CU^{**}(S\cup(\omega_1-S_\eta))$.''

{\sl Then $\langle P_\eta\,\colon\eta\leq\kappa\rangle$ is strictly $(T,S)$-$\#$-preserving.}

\medskip

Comment: Note that the stationary sets which are collapsed in case (2) are not in 
the ground model.  Indeed, we have that $P_\alpha$ is proper for $\alpha\leq\kappa$.

Proof:  We prove by induction on $\kappa$ that whenever $\lambda$ is a sufficiently large regular cardinal and
$N$ is a countable elementary substructure of $H_\lambda$ and $\{\langle( {\cal I}_\eta,S_\eta)\,\colon
\allowbreak\eta<\kappa\rangle,P_\kappa,T,S\}\in N$ and $\alpha\in\kappa\cap N$ and
$q$ is $(N,P_\alpha,S,T)$-$\#$-preserving and $q\forces``\dot p\in\dot P_{\alpha,\kappa}\cap N[G_{P_\alpha}]$''
then there is $r\in P_\kappa$ such that $r\restr\alpha=q$ and $q\forces``r\restr[\alpha,\kappa)\leq\dot p$'' and
$\supt(r)\subseteq\alpha\cup N$ and $r$ is $(N,P_\kappa,S,T)$-$\#$-preserving.  This differs 
from the definition of
$(T,S)$-$\#$-preserving insofar as we assume that $N$ must contain a certain additional
 parameter (namely, $\langle({\cal I}_\eta,S_\eta)\,\colon\allowbreak\eta<\kappa\rangle$), 
but by Lemma 5.3 this is immaterial.

In line with the induction on $\kappa$, we assume that $P_\beta$ is strictly
$(T,S)$-$\#$-preserving for every $\beta<\kappa$.  Assume $\lambda$ and $N$ are given.

Case 1 (successor step).  Suppose $\kappa=\eta+1$.

Necessarily $\eta\in N$, so we may assume that $\eta=\alpha$.  If ${\bf 1}\forces_{P_\eta}
``\dot Q_\eta$ is $(T,S)$-$\#$-preserving,'' then we are done by Lemma 3.5.  Otherwise, by 
Lemmas 3.4 and 3.11 we have that 
$q\forces``\omega_1\cap N[G_{P_\eta}]=\omega_1\cap N\notin S_\eta$,'' so by
Lemmas 3.5 and 3.9 we are again done.

Case 2 (limit step). Suppose $\kappa$ is a limit ordinal.

Let $\kappa'={\rm sup}(\kappa\cap N)$, and let $\langle\alpha_i\,\colon\allowbreak 
i\in\omega\rangle$ be a strictly increasing sequence of ordinals from $\kappa\cap N$ 
cofinal in $\kappa'$ such that $\alpha_0=\alpha$.  Let $\langle\sigma_i\,\colon i\in
\omega\rangle$ list
all $P_\kappa$-names $\sigma$ in $N$ such that ${\bf 1}\forces_{P_\kappa}``\sigma$ is
an ordinal.''  Let $\langle(x_n,A_n)\,\colon\allowbreak n\in\omega\rangle$
list all pairs $\langle x,A\rangle$ 
such that $A\in N$ is a $P_\kappa$-name for a subset of $T$  and
$x$ is $(T,N)$-$\#$-generic.  
Build a sequence $\langle(q_n,\dot p_n)\,\colon\allowbreak n\in\omega\rangle$ such that
$q_0=q$ and $\dot p_0=\dot p$ and for every $m\in\omega$ we have

(1) $q_m$ is $(N,P_{\alpha_m},S,T)$-$\#$-preserving and ${\rm supt}(q_m)\subseteq
\alpha\cup N$,

(2) $q_m\forces``q_{m+1}\restr[\alpha_m,\alpha_{m+1})\leq\dot p_m\restr\alpha_{m+1}$'' 
and $q_{m+1}$ is
$(N,\allowbreak P_{\alpha_{m+1}},S,T)$-$\#$-preserving,

(3) $q_{m+1}\restr\alpha_m=q_m$,

(4) $q_{m+1}\forces``\dot p_{m+1}\leq \dot p_m\restr[\alpha_{m+1},\kappa)$ and
$\dot p_{m+1}\in N[G_{P_{\alpha_{m+1}}}]$ and $\dot p_{m+1}$ decides the value of
$\sigma_m$,''

(5) 
$q_{m+1}\forces``\dot p_{m+1}\forces`(\exists y<x_m)\allowbreak
(y\in A_m$ or $(\forall z\geq y)\allowbreak(z\notin A_m))$'\thinspace''

This is possible by Lemma 3.10 and the induction hypothesis.

Take $r\in P_\kappa$ such that ${\rm supt}(r)\subseteq\alpha\cup N$ and for every
$m\in\omega$ we have that $r\restr\alpha_m=q_m$.  
This concludes the induction, and thereby establishes the Lemma.

\section{The models}

Shelah (item (2) below) and, later, Schlindwein, have constructed models of each of the following:

(1) every Aronszajn tree is $S$-$*$-special ($S$ an arbitrary stationary set
that is in the ground model) and some Aronszajn tree
$T$ is not $S'$-$*$-special whenever $S'-S$ is stationary, and CH holds [JSL],

(2) every Aronszajn tree is $S$-$*$-special 
($S$ an arbitrary stationary set that is in the ground model) and
some Aronszajn tree $T$ is not $S'$-$*$-special whenever $S'-S$ is stationary, and CH fails ([Sh], or use
[JSL] with $CU^{**}(S\cup(\omega_1-S_\alpha))$ in place of
$CU(S\cup(\omega_1-S_\alpha))$),

(3) Suslin's hypothesis plus some Aronszajn tree $T$
 has no stationary antichain, plus CH fails [APAL],

(4) Suslin's hypothesis plus some Aronszajn tree $T$
has no stationary antichain, plus CH holds [APAL2],

(5) every $\omega_1$-tree is $S$-$*$-special for $S$ an arbitrary
 stationary set
that is in the ground model (in particular, Kurepa's hypothesis fails)
and some $\omega_1$-tree
$T$ is not $S'$-$*$-special whenever $S'-S$ stationary,
 plus CH holds [STACY],

(6) same as (5) but CH fails
 (use [STACY] but with $CU^{**}(S\cup(\omega_1-S_\alpha))$ in place of
  $CU(S\cup(\omega_1-S_\alpha))$).

Models (5) and (6) require an inaccessible cardinal.

We claim that in variants of each of these six models 
(seven, actually, as there are two different constructions cited in item (2)) 
there is no unbounded $E\subseteq\omega_1$ such that
$T$ is $E$-special.
The demonstrations are all entirely similar to each other, 
except in Shelah's construction for item (2), where the
demonstration of $(S,T)$-$\#$-preserving for the forcing of
[PIF XXXXX] is similar to the proof of Lemma 3.9 (second case) above.
Because the changes to previously published material 
are easily explained, we do not give complete proofs for all six
models.

\section{The first two models}

In this section, we show that a variation  of the model from [JSL]
satisfies that there is an Aronszajn tree $T^*$ and a stationary set $S^*$
such that every Aronszajn tree is $S^*$-$*$-special, and whenever
$S'-S^*$ is stationary then $T^*$ is not $S'$-$*$-special, and
for every unbounded $E\subseteq\omega_1$ we have that $T^*$ is not
$E$-special. The model in [JSL] was used to solve the problem of
constructing a model of
ZFC plus CH plus SH plus not every Aronszajn tree is special
(answering a  question posed by Shelah [PIF, XXX]).

Throughout this section, we fix $S^*$ a stationary co-stationary subset of
$\omega_1$ and we fix an Aronszajn tree $T^*$.  In the end, we will use
a sharply Suslin tree in the ground model as $T^*$, so that
the final poset $P_{\omega_2}$ will force for
 every $S'\subseteq\omega_1$ such that
$S'-S^*$ is stationary, then $T^*$ is
not $S'$-$*$-special, and for every unbounded $E\subseteq\omega_1$
we have that $T^*$ is not $E$-special.

For $T$ an Aronszajn tree, we let $T^n$ be the Aronszajn tree consisting of $N$-tuples
of elements of $T$, all of which have the same rank.
For $x\in T$ and $\beta\leq{\rm rk}(x)$ we let
$x\restr\beta$ be the unique $y\leq x$ such that
${\rm rk}(y)=\beta$.
We turn our attention to defining the posets
that will be used as the constituent posets of the iteration.

\proclaim Definition 7.1. We say that $R$ is a finite rectangle iff there is some $n=n(R)\in\omega$
and some sequence $\langle R_i\,\colon\allowbreak i<n\rangle$
such that $R=R_0\,{\sf X}\,R_1\,{\sf X}\cdots{\sf X}\,R_{n-1}$ and
for each $i<n$ we have that $R_i$ is a finite subset of $\omega_1$.

\proclaim Definition 7.2. Suppose $T$ is an Aronszajn tree
and $\gamma<\omega_1$ and $n\in\omega$ and
 ${\overline x}\in T^n_\gamma$ and $f$ is an ordinal-valued function
and $R=R_0\,{\sf X}\,R_1\,{\sf X}\cdots{\sf X}\,R_{n-1}$ is a finite rectangle.
Then we define $\heartsuit(\alpha,{\overline x},f,R)$ to mean
that whenever $\alpha<\beta\leq\gamma$ and $i<n$ and
$x_i\restr\beta\in\dom(f)$ then $f(x_i\restr\beta)\notin R_i$.

\proclaim Definition 7.3.  Suppose $T$ is an Aronszajn tree.  We let $P'(T)$ be
the poset whose universe is\/ {\rm $\{\langle f,S\rangle\,\colon\allowbreak
S$ is a countable set of countable limit ordinals and $f$ is an
 $S$-$*$-specializing  function, and ${\rm cl}(S)\cap S^*\subseteq S\}$.}
 The ordering is given by co-ordinatewise reverse end-extension.

\proclaim Definition 7.4. $\Gamma$ is  a $T$-promise iff there is a
closed unbounded $C=C(\Gamma)\subseteq\omega_1$ and an integer $n=n(\Gamma)$ and an ${\overline x}
={\rm min}(\Gamma)\in \Gamma$ such that
$\Gamma\subseteq T^n_C$ and whenever $\alpha<\beta$
are in $C$ and ${\overline y}\in \Gamma\cap T^n_\alpha$ then there
is an infinite $W\subseteq\Gamma\cap T^n_\beta$ such that for
every ${\overline z}\in W$ we have ${\overline x}\leq
{\overline y}\leq{\overline z}$, and distinct elements of $W$ have
disjoint ranges.  We also require that for every ${\overline z}\in\Gamma$
and every
$\alpha\in {\rm rk}({\overline z})\cap C$ we have
${\overline z}\restr\alpha\in\Gamma$.

The following Fact is proved in [PIF XXXX] and [JSL, Lemma 50].

\proclaim Fact 7.5. Suppose $\Delta\subseteq T^n$ is uncountable and downwards closed and every element of $\Delta$ is
comparable with ${\overline x}\in T^n$. Then there is some $T$-promise
$\Gamma\subseteq\Delta$ such that ${\rm min}(\Gamma)={\overline x}$.

\proclaim Definition 7.6.  Suppose $\langle f,S\rangle\in P'(T)$ and $\Gamma$ is a
$T$-promise.  We say that $\langle f,S\rangle$ fulfills $\Gamma$ iff
$S-{\rm rk}({\rm min}(\Gamma))\subseteq C(\Gamma)$ and whenever
$\beta\in C(\Gamma)$ and $\alpha\in C(\Gamma)\cap S\cap\beta$
and ${\overline y}\in\Gamma\cap T^n_\beta$ and $R$ is a
finite rectangle with $n(\Gamma)=n(R)$ then
there is an infinite $W\subseteq\Gamma\cap T^n_\beta$ such that
distinct elements of $W$ have disjoint ranges and
for every ${\overline w}\in W$ we have ${\overline y}\leq
{\overline w}$ and $\heartsuit(\alpha,{\overline w},f,R)$.

Note that in  Definition 7.6 we do not assume that
$\beta\in S$.

Contrast the following Definition with [JSL, Definition 52]. The difference is that
in [JSL], it is required that $\Psi$ be countable.

\proclaim Definition 7.7. $P(T)$ is the poset whose universe consists of
triples $\langle f,S,\Psi\rangle$ such that
$\langle f,S\rangle\in P'(T)$ and $\Psi $ is a set
of $T$-promises that $\langle f,S\rangle$ fulfills such that
for every $\alpha<\omega_1$ we have that $\{\Gamma\in\Psi\,\colon\allowbreak
{\rm rk}({\rm min}(\Gamma))<\alpha\}$ is countable.
The ordering is given by $\langle f',S',\Psi'\rangle\leq
\langle f,S,\Psi\rangle$ iff $f\subseteq f'$ and
$S'\cap({\rm sup}(S)+1)=S$ and $\Psi\subseteq\Psi'$.

For $p\in P(T)$ we will use $f_p$, $S_p$, and $\Psi_p$ to denote the components of $p$, and we set
${\rm ht}(p)$ to equal ${\rm sup}(S_p)$.  We set $C(\Psi_p)=
\{\gamma<\omega_1\,\colon(\forall\Gamma\in\Psi_p)\allowbreak
({\rm rk}({\rm min}(\Gamma))>\gamma$ or
$\gamma\in C(\Gamma))\}$. Note that $C(\Psi_p)$ is closed and
unbounded.

The following Lemmas correspond to [JSL, Lemmas 53 through 56].  Although the poset $P(T)$ referred to in [JSL] differs
from $P(T)$ in that [JSL] required $\Psi_p$ is countable for every $p\in P(T)$,
the difference is immaterial to the proof of these Lemmas.

\proclaim Lemma 7.8.  Suppose $p\in P(T)$ and ${\rm ht}(p)=\alpha<\beta\in
C(\Psi_p)$.  Suppose
$R$ is a finite rectangle and ${\overline z}\in T^n_\beta$. Then there is
$q\leq p$ such that ${\rm ht}(q)=\beta$ and
$\heartsuit(\alpha,{\overline z},f_q,R)$.

Proof: In the proof of [JSL, Lemma 53] simply replace the clause
$``\Gamma\in\Psi_p$'' with $``\Gamma\in\Psi_p$ and
${\rm rk}({\rm min}(\Gamma))\leq\beta$.'' Besides that change, the proof is
unchanged.

The following Lemma is [JSL, Lemma 54] but, as usual, for a slightly different poset.
The same Lemma basically appears in [PIF, chapter V] for a different poset.

\proclaim Lemma 7.9. Suppose that $\lambda$ is a sufficiently large regular cardinal and $N$ is
a countable elementary substructure of $H_\lambda$ and
$P(T)\in N$ and $p\in P(T)\cap N$ and $D\in N$ is dense in $P(T)$
and $\delta=\omega_1\cap N$ and $R$ is a finite rectangle and
$n=n(R)$ and ${\overline x}\in T^n_\delta$. Then there is
$q\leq p$ such that $q\in D\cap N$ and $\heartsuit({\rm ht}(p),{\overline x},f_q,R)$.

Proof: The proof of [JSL, Lemma 54] carries over verbatim.

\proclaim Lemma 7.10. Suppose that $\lambda$ is a sufficiently
large regular cardinal
and $N$ is a countable elementary substructure of $H_\lambda$ containing $P(T)$.
Let $\delta=\omega_1\cap N$, and suppose $p\in P(T)\cap N$.
Suppose also that ${\overline x}\in T^n_\delta$ and $R$ is a
finite rectangle and $n=n(R)$.  Then there is $q\leq p$
such that $\heartsuit({\rm ht}(p),
{\overline x},f_q,R)$ and ${\rm ht}(q)=\delta$ and for every
open dense $D\in N$ there is some $r\in D\cap N$ such that
$q\leq r$.
In particular, $P(T)$ is proper and does not add reals.

Proof: The proof of [JSL, Lemma 55] carries over verbatim.

For the definition of $(S^*,\omega_2)$-p.i.c., see [JSL,Definition 39].

\proclaim Lemma 7.11. $P(T)$ has $(S^*,\omega_2)$-p.i.c.

Proof: The proof of [JSL, Lemma 57] carries over verbatim.  However, please note the error in
[JSL, Lemma 43] concerning the preservation of
$(S^*,\omega_2)$-p.i.c.  The corrected statement of that Lemma is
as follows:

\proclaim Lemma 7.12. Suppose $\langle P_\eta\,\colon\eta\leq\kappa\rangle$
is a countable support forcing iteration based on
$\langle Q_\eta\,\colon\allowbreak\eta<\kappa\rangle$.  Suppose that for every
$\eta<\kappa$ we have that
${\bf 1}\forces_{P_\eta}``Q_\eta$ {\rm has $(S,\omega_2)$-p.i.c.''}
Then if $\kappa<\omega_2$ we have that $P_\kappa$ has
$(S,\omega_2)$-p.i.c., and if $\kappa\leq\omega_2$ then
$P_\kappa$ has $\omega_2$-c.c.

Proof: [JSL, Lemma 43] neglects the restriction on the
 length of the iteration, and the proof given there is incorrect.
The needed correction to the proof is to be found in
the proof of [APAL, Lemma XXXX].

 \proclaim Lemma 7.13. $ {\bf 1}\forces_{P(T)}``T$ {\rm is
 $S^*$-$*$-special.''}

 Proof: The proof of [JSL, Lemma 58] carries over verbatim.

The following is [PIF, remark on page XXX], [JSL, Lemma 8].

\proclaim Lemma 7.14. Suppose $\langle P_\eta\,\colon\allowbreak\eta\leq\kappa\rangle$ is as in
Lemma XXX, and suppose\/ {\rm
${\bf 1}\forces_{P_1}``T$ is $S$-$*$-special.''} Then\/ {\rm
 ${\bf 1}\forces_{P_\kappa}
``T$ is Aronszajn.''}

Proof: See [JSL, Lemma 8].

We now turn to the task of showing that $P(T)$ is $(T^*,S^*)$-$\#$-preserving,
and hence by Lemmas 4.3, 4.4, and 4.5, and Theorem 5.4 we may construct
the iteration so that in $V[G_{P_{\omega_2}}]$ we have:

(1) $T^*$ is Aronszajn,

(2) for every Aronszajn tree $T$ we have $T$ is $S^*$-$*$-special,

(3) for all $S'\subseteq\omega_1$ such that $S'-S^*$ is stationary, we have
that $T^*$ is not $S'$-$*$-special,

(4) for all unbounded $E\subseteq\omega_1$ we have that
$T^*$ is not $E$-special.

Simply take $T^*$ to be sharply Suslin in the ground model, and choose $Q_0$
to be $P(T^*)$ so that the hypotheses of Lemmas 4.5 and
7.14 are satisfied. In order to ensure (3), use posets of the form
$CU(S^*\cup(\omega_1-S_\eta))$ or
$CU^*(S^*\cup(\omega_1-S_\eta))$ or
$CU^{**}(S^*\cup(\omega_1-S_\eta))$, where
$S_\eta=\{{\rm rk}(x)\,\colon\allowbreak x\in{\cal I}_\eta\}$ where
${\cal I}_\eta$ is (a name for) an antichain of $T^*$.

\proclaim Lemma 7.15 (analogue of [JSL, Lemma 59]).
  Suppose $P=P(T)$ and $\lambda$ is a sufficiently large
  regular cardinal and $N$ is a countable
  elementary substructure of $H_\lambda$ containing
  $\{P,S^*,T^*\}$. Let $\delta=\omega_1\cap N$.
  Suppose $m\in\omega$ and
  ${\overline x}\in T^m_\delta$ and $p\in P\cap N$
and $R\in N$ is a finite rectangle and ${\overline z}\leq{\overline x}$
and ${\rm rk}({\overline z})={\rm ht}(p)$ and $A\in N$ is a
 $P$-name for a subset of
$T^*$ and $x$ is $(T^*,N)$-$\#$-generic. Then there is $y<x$
and $q\leq p$ such that
$q\in N$ and
$\heartsuit({\rm ht}(p),{\overline x},f_q,R)$ and either\/ {\rm
$q\forces``y\in A$''} or\/ {\rm
$q\forces
``(\forall z\geq y)\allowbreak(z\notin A)$.''}

Proof: Let $A^*=\{y\in T^*\,\colon\allowbreak$ for every $T$-promise
$\Gamma$ such that ${\rm rk}({\rm min}(\Gamma))\geq{\rm max}(
{\rm ht}(p),\allowbreak{\rm rk}(y))$ we have
$\langle f_p,S_p,\Psi_p\cup\{\Gamma\}\rangle\notforces``y\notin A$''$\}$.
Notice $A^*\in N$, and therefore we may fix
$y<x$ such that either $y\in A^*$ or $(\forall z\geq y)\allowbreak
(z\notin A^*)$.

Case 1.  $y\in A^*$.

Suppose there is no $q\leq p$ such that
$q\forces``y\in A$'' and ${\rm ht}(q)<\delta$ and
$\heartsuit({\rm ht}(p),{\overline x},f_q,R)$.
Fix $\alpha\in C(\Psi_p)$ such that $\alpha\geq
{\rm max}({\rm ht}(p),\allowbreak
{\rm rk}(y))$.
Let $\Delta=\{{\overline w}\in T^n\,\colon\allowbreak$there is
no $q\leq p$ such that $q\forces``y\in A$'' and
${\rm ht}(q)<{\rm rk}({\overline w})$ and $\heartsuit({\rm ht}(p),
{\overline w},f_q,R)$ and
${\overline w}$ is comparable with
${\overline x}\restr\alpha\}$. Notice $\Delta\in N$. We have that every ${\overline w}\leq
{\overline x}$ is in $\Delta$. Hence

\medskip

{\centerline{$N\models``\Delta$ is uncountable.''}}

\medskip

We also have that $\Delta$ is downwards closed.
Hence by Fact 7.5 we may take $\Gamma\subseteq\Delta$ such that
$\Gamma$ is a $T$-promise and ${\rm min}(\Gamma)\leq{\overline x}$
and ${\rm rk}({\rm min}(\Gamma))\geq{\rm max}({\rm ht}(p),\allowbreak
{\rm rk}(y))$.

Because $y\in A^*$ we have that
$\langle f_p,S_p,\Psi_p\cup\{\Gamma\}\rangle\notforces``y\notin A$.''
Therefore we may take $r\leq\langle f_p,S_p,\Psi_p\cup\{\Gamma\}\rangle$
such that $r\forces``y\in A$.''   Because $\langle f_r,S_r\rangle$
fulfills $\Gamma$, we may take ${\overline w}\in\Gamma$ such that
${\rm rk}({\overline w})>{\rm ht}(r)$ and
$\heartsuit({\rm ht}(p),{\overline w},f_r,R)$.
Because ${\overline w}\in\Delta$ there is no $q\leq p$ such that
$q\forces``y\in A$'' and ${\rm ht}(q)<{\rm rk}({\overline w})$ and
$\heartsuit({\rm ht}(p),{\overline w},f_q,R)$.  But $r$ witnesses the opposite.
This contradiction shows that 
if Case 1 holds then there is  $q\leq p$ such that
$q\forces``y\in A$'' and ${\rm ht}(q)<\delta$ and
$\heartsuit({\rm ht}(p),{\overline x},f_q,R)$.
Let $\beta={\rm ht}(q)$.
We have $(\exists q\leq p)\allowbreak({\rm ht}(q)=\beta$ and
$\heartsuit(\alpha,{\overline x}\restr\beta,\allowbreak
f_q,R)$ and $q\forces``y\in A$'').  Because
${\overline x}\restr\beta\in N$ we have
 $(\exists q\leq p)\allowbreak({\rm ht}(q)=\beta$ and $q\in N$ and
$\heartsuit(\alpha,{\overline x}\restr\beta,\allowbreak
f_q,R)$ and $q\forces``y\in A$'').  Hence the
 conclusion of the Lemma holds in
Case 1.

Case 2: Otherwise.

We have that $(\forall z\geq y)\allowbreak
(z\notin A^*)$.  Hence
 $(\forall z\geq y)\allowbreak(\exists\Gamma(z))\allowbreak
(\Gamma(z)$ is a $T$-promise and ${\rm rk}({\rm min}(\Gamma))\geq
{\rm max}({\rm rk}(z),{\rm ht}(p))$ and
$\langle f_p,S_p,\Psi_p\cup\{\Gamma(z)\}\rangle\forces``z\notin A$'').
We may assume that the function mapping $z$ to $\Gamma(z)$ is an
element of $N$.
Let $q=\langle f_p,S_p,\Psi_p\cup\{\Gamma(z)\,\colon\allowbreak
z\geq y\}\rangle$.  We have
$q\in P(T)\cap N$ and
$q\forces``(\forall z\geq y)\allowbreak(z\notin A$).''
Hence the conclusion of the Lemma holds in Case 2.

The Lemma is established.

\proclaim Theorem 7.16.  $P(T)$ is $(T^*,S^*)$-$\#$-preserving.

Proof: Suppose $\lambda$ is a sufficiently large regular cardinal and $N$
is a countable elementary substructure of $H_\lambda$ containing $
\{P(T),T^*,S^*\}$. Let $\delta=\omega_1\cap N$. If $\delta\in S^*$
then we are done by Lemma 7.10, so assume otherwise.
Similarly, if there are no $x\in T^*_\delta$ such that
$x$ is $(N,T^*)$-$\#$-generic, we are done, so
assume otherwise. Suppose $p\in P(T)\cap N$. We build $q\leq p$
such that $q$ is $(N,P(T),S^*,T^*)$-$\#$-preserving. Let
$\langle D_m\,\colon\allowbreak m\in\omega\rangle$ list the dense open subsets of
$P(T)$ that are in $N$. Let
$\langle(\Gamma_m,R'_m,{\overline z}_m)\,\colon\allowbreak m\in\omega\rangle$
list all triples $\langle \Gamma, R', {\overline z}\rangle$ such that
$\Gamma\in N$ is a $T$-promise and $R'\in N$ is a finite rectangle and
${\overline z}\in \Gamma\cap N$ and $n(R')=n(\Gamma)$, with infinitely
many repetitions.  Let $\langle (x_m,A_m)\,\colon\allowbreak
m\in\omega\rangle$ list all pairs $(x,A)$ such that
$x$ is $(T^*,N)$-$\#$-generic and $A$ is a $P(T)$-name in $N$ for
a subset of $T^*$.

Build by recursion $\langle(F_m,q_m,p_m,{\overline w}_m)\,\colon
\allowbreak m\in\omega\rangle$ such that $F_0=\emptyset$
and $p_0=p$ and  each of the following holds:

(1)  $F_m$ maps a finite subset of $T_\delta$ into the set of
finite subsets of $\delta$,

(2) $q_m\in P(T)\cap N$ and $q_m\leq p_m$ and $(\forall w\in
{\rm dom}(F_m))\allowbreak(\heartsuit({\rm ht}(p_m),w,f_{q_m},
F_m(w)))$ and
for some $y<x_m$ we have either
$q_m\forces``y\in A_m$'' or $q_m\forces``(\forall z\geq y)\allowbreak
(z\notin A_m)$,''

(3) $p_{m+1}\in D_m\cap N$ and $p_{m+1}\leq q_m$ and
$(\forall w\in{\rm dom}(F_m))\allowbreak(\heartsuit
({\rm ht}(q_m),w,f_{p_{m+1}},F_m(w)))$,

(4) if $\Gamma_m\in\Psi_{p_{m+1}}$ and ${\rm range}({\overline z}_m)
\subseteq {\rm dom}(f_{p_{m+1}})$ then
${\overline w}_m\in\Gamma_m\cap T^{n(\Gamma_m)}_\delta$ and
${\overline z}_m\leq{\overline w}_m$ and for all $i<n(\Gamma_m)$ we have
$\heartsuit({\rm rk}({\overline z}_m),{\overline w}_m(i),
f_{p_{m+1}},(R'_m)_i)$ and ${\rm range}({\overline w}_m)$ is disjoint
from ${\rm dom}(F_m)$; otherwise, ${\overline w}_m=\emptyset$,

(5) ${\rm dom}(F_{m+1})={\rm dom}(F_m)\cup{\rm range}({\overline w}_m)$,

(6) for every $w\in{\rm dom}(F_m)$ we have $F_{m+1}(w)\supseteq
F_m(w)$,

(7) for all $i\in{\rm dom}({\overline w}_m)$ we have $F_{m+1}(
{\overline w}_m(i))\supseteq (R'_m)_i$.

The construction can be carried out using Lemma 7.15 to choose $q_m$ as in (2) and
Lemma 7.9 to choose $p_{m+1}$ as in (3).

Let $q=\langle\bigcup\{f_{p_m}\,\colon\allowbreak m\in\omega\},\allowbreak
\bigcup\{S_{p_m}\,\colon\allowbreak m\in\omega\},\allowbreak
\bigcup\{\Psi_{p_m}\,\colon\allowbreak m\in\omega\}\rangle$. It is easy to see
that $q$ is as required.

The theorem is established.

We turn to the problem of showing that the  forcing iteration under
consideration does not add reals.

For the definition of ``$(S,<\omega_1)$-proper,'' see
 [JSL, Definition 26].

\proclaim Lemma 7.17.  $P(T)$ is $(\omega_1,<\omega_1)$-proper,
and hence $P(T)$ is $(S,<\!\omega_1)$-proper for any stationary
$S\subseteq\omega_1$.

The proof of [JSL, Lemma 56] carries over verbatim, but there is a small error:
namely, in Case 1 it is implicitly assumed that $\gamma\ne 0$.  Fortunately,
the case that $\gamma=0$ is easily handled by
simply setting $r^*=p$.

The following is [JSL, Definition 27], based on [PIF, Chapter V].

\proclaim Definition 7.18. Suppose $\langle P_\eta\,\colon\allowbreak\eta\leq\kappa\rangle$ is
a countable support forcing iteration. We say that
the iteration is strictly $(S,<\!\omega_1)$-proper iff {\bf whenever}
$\rho<\omega_1$ and $\lambda$ is a sufficiently large
regular cardinal and $\langle N_i\,\colon\allowbreak i\leq\rho\rangle$ is
a continuous tower of countable elementary substructures of $H_\lambda$ and
$P_\kappa\in N_0$ and  for every $i<\rho$ we have
$\langle N_j\,\colon\allowbreak j\leq i\rangle\in N_{i+1}$ and
for every $i\leq\rho$ we have
$\omega_1\cap N_i\in S$ and $i\in N_i$ and $\eta\in \kappa\cap N_0$ and
$p\in P_\eta$ and for every $i\leq\rho$ we have that $p$ is
$(N_i,P_\eta)$-generic, and $p\forces``q\in P_{\eta,\kappa}\cap
N_0[G_{P_\eta}]$,'' {\bf then} there is $r\in P_\kappa$ such that
$r\restr\eta=p$ and $p\forces``r\restr[\eta,\kappa)\leq q$'' and
for every $i\leq\rho$ we have that $r$ is $(N_i,P_\kappa)$-generic and
${\rm supt}(r)\subseteq\eta\cup N_\rho$.

\proclaim Definition 7.19. Suppose $\lambda$ is
large for $P$ and $M$ is a countable elementary
substructure of\/ $H_\lambda$ and
$P\in M$ and $p\in P\cap M$.  We set
$\Gen(M,P,p)$ equal to the set of all $G\subseteq P\cap M$
which satisfy all of the following:

\noindent(1) $G$ is $M$-generic, i.e., whenever $D\in M$ is a dense
open subset of $P$ then $G\cap D\ne \emptyset$

\noindent(2) $G$ is directed, i.e., $(\forall q_1\in G)\allowbreak
(\forall q_2\in G)\allowbreak(\exists r\in G)\allowbreak
(r\leq q_1$ and $r\leq q_2)$

\noindent(3) $p\in G$ 

\proclaim Definition 7.20.  Suppose that $S$ is stationary and
$P$ is $S$-proper not adding reals and
suppose $Q$ is a $P$-name for a poset.
We say $Q$ is $S$-complete for $P$ iff whenever
$\lambda$ is a sufficiently large regular cardinal
and\/ $M$ and $N$ are countable elementary
substructures of\/ $H_\lambda$ and\/ $P*Q\in M\in N$ and
$\omega_1\cap M\in S$ and $\omega_1\cap N\in S$
and $q\in M$ is a $P$-name for an element of $Q$
and $G\in\Gen(M,P,{\bf 1})\cap N$, then there is
$G'\in{\rm Gen}(M,P*Q,{\bf 1})$
such that $\{p_1\in P\,\colon\allowbreak
(\exists r)\allowbreak((p_1,r)\in G')\}=G$ and
$({\bf 1},q)\in G'$ and whenever $p$ is a lower
bound for $G$ and $p$ is $N$-generic, then there
is $p'\leq p$ and a $P$-name $s$ for an element of $Q$ such that
$(p',s)$ is a lower bound for $G'$.

\proclaim Lemma 7.21. Suppose $S^*\subseteq S$.
Then $CU(S)$ and $CU^*(S)$ are $(S^*,<\omega_1)$-proper.
Furthermore, for every $P$ such that $P$ is $S^*$-proper not
adding reals, if\/ {\rm ${\bf 1}\forces_P``Q=CU(S)$ or $Q=CU^*(S)$ and
$S^*\subseteq S$''} then $Q$ is $S^*$-complete for $P$.

Proof: See [JSL, Lemmas 37 and 38].

\proclaim Lemma 7.22. Suppose $S\subseteq\omega_1$ is stationary
and $\langle P_\eta\,\colon\allowbreak\eta\leq\kappa\rangle$
is a countable support iteration based on
$\langle Q_\eta\,\colon\allowbreak\eta<\kappa\rangle$
and for every $\eta<\kappa$ we have that
$Q_\eta$ is $S$-complete for $P_\eta$, and suppose
also that $P_\kappa$ is strictly $(S,<\!\omega_1)$-proper. Then
$P_\kappa$ does not add reals.

Proof: See [JSL, Theorem 36].

\proclaim Theorem 7.23.  If ZFC is consistent, then so is ZFC plus
there is a stationary co-stationary set $S^*$ such that every Aronszajn tree
is $S^*$-$*$-special plus there is an Aronszajn tree $T^*$ such that
$T^*$ is not $S$-$*$-special whenever
$S-S^*$ is stationary, and for every unbounded
$E\subseteq\omega_1$ we have that $T^*$ is not $E$-special.
Furthermore, we may either have CH hold or CH fail in the model.

Proof: This is [JSL, Theorem 45] with three changes.  The first change is
that we allow forcings of the form $CU^{*}(S^*\cup(\omega_1-S_\eta))$ and
 $CU^{**}(S^*\cup(\omega_1-S_\eta))$. Naturally,
 by using  $CU^{**}(S^*\cup(\omega_1-S_\eta))$ we will
 not have CH in the final model. The second change is that
we start with a sharply Suslin tree and have the property of
$(S^*,T^*)$-$\#$-preserving in order to assure that  for every
unbounded $E\subseteq\omega_1$ we have\/ {\rm
${\bf 1}\forces_{P_{\omega_2}}``T^*$ is not $E$-special.''}
The third change is that we use the version of $P(T)$ in which
$\Psi$ is not required to be countable, but only that for each
$\alpha$ we have that
$\{\Gamma\in\Psi\,\colon\allowbreak{\rm rk}({\rm min}(\Gamma))\leq
\alpha\}$ is countable.

\section{The ``no stationary antichains'' models}

In [APAL] a model of ZFC plus SH plus some Aronszajn tree has no sttionary antichain
is constructed.  The Suslin trees are killed more gently than in
[PIF, Chapter IX.4] and [PIF] and [STACY], because there is an Aronszajn tree
$T^*$ such that for every stationary $S\subseteq\omega_1$ we have that
$T^*$ is not $S$-$*$-special.  Here we show how to kill Suslin trees even
more gently; in the final model, Suslin's hypothesis holds and there is
an Aronszajn tree $T^*$ such that for every stationary $S$ and every unbounded $E$
we have that $T$ is neither $S$-$*$-special nor $E$-special.

First we recall the differences between the construction of [JSL] and [APAL].

The first difference is the poset that is used.  Before
exhibiting the poset from [APAL], we give some definitions.

\proclaim Definition 8.1.  Suppose $T$ is an Aronszajn tree and
$f$ is a monotonically non-decreasing
function from $\bigcup\{T_\beta\,\colon\allowbreak\beta\leq
\alpha\}$ into $\{0,1\}$.  We set ${\rm ht}(f)$ equal to $\alpha$.
Given ${\overline z}\in T^n$ and $\rho<\omega_1$, we say
$\heartsuit(\rho,f,{\overline z})$ iff either
 $\rho\geq{\rm rk}({\overline z})$ or for all $i<n$ and all $t\leq
{\overline z}(i)$ such that $t\in{\rm dom}(f)$
we have $f(t)=f({\overline z}(i)\restr\rho)$.
Given a $T$-promise $\Gamma$, we say that $f$ fulfills $\Gamma$
iff whenever $\beta<\gamma$ are in
$C(\Gamma)$ and $\beta<{\rm ht}(f)$ and ${\overline w}\in\Gamma\cap
T^{n(\Gamma)}_\beta$, then there is an infinite $W\subseteq
\Gamma\cap T^n_\gamma$
such that distinct elements of $W$ have disjoint ranges and
for every ${\overline w}\in W$ we have
$\heartsuit(\beta, f,{\overline w})$.

In [APAL] we may view the
poset $P(T)$ as $\{\langle f,\Psi\rangle\,\colon\allowbreak$ for some
$\alpha<\omega_1$ we have that $f$ is a monotonically non-decreasing
function from $\bigcup\{T_\beta\,\colon\allowbreak\beta\leq\alpha\}$ into
$\{0,1\}$ and $\Psi$ is a countable set of promises that $f$ fulfills$\}$.

We have ${\bf 1}\forces_{P(T)}``T$ is not Suslin because
$\{t\in T\,\colon\allowbreak t$ has an immediate predecessor $t'$ such that
$f(t')=0$ and $f(t)=1$, where $f=\bigcup\{f'\,\colon\allowbreak
(\exists\Psi)\allowbreak(\langle f,\Psi\rangle\in G_{P(T)})$.''

The second change is that we must ensure that the tree $T^*$ must remain Aronszajn
in $V[G_{P_{\omega_2}}]$. Recall that in [JSL] (and in [PIF, Chapter IX.4])
this was accomplished by $S^*$-$*$-specializing $T^*$ at the first step of the
iteration, so that it could not become non-Aronszajn in any extension
in which $\omega_1$ is not collapsed.  This strategy is not available in the
construction in [APAL].  Instead, the fact that $T^*$ remains Aronszajn is ensured by
showing that the iteration satisfies a preservation property that is
more stringent than the property $(T^*,S^*)$-preserving used in
[PIF, Section IX.4] and [JSL].  See [APAL, Definition XXX] for the
definition of this property, [APAL, Lemma XXX] for the fact that
the property ensures that $T^*$ remains Aronszajn, and various
Lemmas in [APAL] for the fact that the property is
preserved under the appropriate iterations and
that the porperty is satisfied by the constituent posets of the iteration.

These two changes were enough to carry out the construction of [APAL], but left the question of
whether we could arrange for CH to hold in the final model.
The difficulty was resolved in [APAL2] by using the following strategy.

\section{Doing without Kurepa trees}

In this section, we strengthen
the conclusion of Theorem 7.23 by requiring that every $\omega_1$-tree
is $S^*$-$*$-special. Therefore we have Kurepa's hypothesis  holds in the
model.  Naturally, this requires  that the hypothesis be strengthened from the
consistency of ZFC to the consistency of ZFC plus there exists an inaccessible
cardinal.  We use (essentially) the forcing from [STACY].

We repeat the main Definitions from [STACY].
We fix $T$ to be an $\omega_1$-tree and
$T^*$  an Aronszajn tree and $B$ equal to the set of uncountable branches
of $T$  and $S^*$ a stationary co-stationary subset of
$\omega_1$.
Fix $\kappa$ a sufficiently large regular cardinal.
For $n\in\omega$ we set $T^n$ equal to $\{{\overline w}\,\colon\allowbreak
{\overline w}$ is a function with domain $n$ and there is some $\alpha<\omega_1$ such that
$(\forall i<n)\allowbreak({\overline w}(i)\in T_\alpha)$. Notice that
for every $n$ we have that $T^n$ is an $\omega_1$-tree.

\proclaim Definition 9.1.  $\Gamma$ is a promise iff
there is $n=n(\Gamma)\in\omega$ and $C=C(\Gamma)\subseteq
\omega_1$ closed unbounded and ${\overline x}={\rm min}(\Gamma)\in
 T^{n(\Gamma)}$
 and $G=G(\Gamma)\subseteq n$
and $\langle b_i(\Gamma)\,\colon\allowbreak i\in G\rangle$ a sequence of
elements of $B$ such that ${\overline x}\in\Gamma\subseteq T^n_C$ and
$(\forall{\overline y}\in\Gamma)\allowbreak({\overline x}\leq{\overline y})$
and for all $\alpha<\beta$ both in $C$ and every
${\overline y}\in\Gamma\cap T^{n(\Gamma)}_\alpha$ then there
is $W\subseteq\Gamma\cap T^{n(\Gamma)}_\beta$ such that
$(\forall{\overline w}\in W)\allowbreak({\overline y}<{\overline w}$ and
$(\forall{\overline w}'\in W)\allowbreak($either ${\overline w}'={\overline w}$
or $\{{\overline w}'(i)\,\colon i\in n(\Gamma)-G\}$ is disjoint from
$\{{\overline w}(i)\,\colon\allowbreak i\in n(\Gamma)-G\}$, and for all
${\overline y}\in\Gamma$ and $i\in G$ we have ${\overline y}(i)\in b_i(\Gamma)$,
and $W$ is infinite unless $G=n$.

Notwithsatnding the fact that we have redefined the notion of ``promise,'' we keep the same definition of ``finite rectangle''
(Definition 7.1).

\proclaim Definition 9.2.  Suppose $n\in\omega$ and ${\overline w}\in T^n$
and $R$ is a finite rectangle and $n(R)=n$ and
$f$ is a function from a subset of $T$ into
$\omega_1$.
Then we say $\heartsuit(\alpha,{\overline w},f,R)$ iff
$(\forall i<n)\allowbreak(\forall y\leq{\overline w}(i))\allowbreak
($if ${\rm rk}(y)>\alpha$ and $y\in{\rm dom}(f)$
then $f(y)\notin R(i))$.
For $b$ an uncountable branch of $T$ we say
$\heartsuit(\alpha,b,f,R)$ iff $n(R)=1$ and
$(\forall x\in b)\allowbreak(\heartsuit(\alpha,x,f,R))$.

\proclaim Definition 9.3.  Suppose $S$ is a bounded subset of $\omega_1$ and
$f$ is a function that $S$-$*$-specializes $T$ and $n\in\omega$ and
$\Gamma$ is a promise and $n=n(\Gamma)$.
We say that $\langle f, S\rangle$ fulfills $\Gamma$ iff
$S-{\rm rk}({\rm min}(\Gamma))\subseteq
C(\Gamma)$ and for every
$\alpha<\beta$ both in $C(\Gamma)$ and every finite rectangle $R$ with
${\rm dom}(R)=n$ and every
${\overline y}\in\Gamma\cap T^{n}_\alpha$ there is 
$W\subseteq\Gamma\cap T^n_\beta$ such that either $G(\Gamma)=n$ or $W$ is infinite, 
and such that
for every ${\overline w}$ and ${\overline w}'$ distinct elements of $W$
we have that $\{{\overline w}(i)\,\colon\allowbreak  i\in n-G(\Gamma)\}\cap
\{{\overline w}'(i)\,\colon\allowbreak  i\in n-G(\Gamma)\}
=\emptyset$ and, for every ${\overline w}\in W$, we have
${\overline y}\leq{\overline w}$ and 
$\heartsuit(\alpha,{\overline w},f,R))$.

\proclaim Definition 9.4. We set $P=P(T,\kappa)$ equal to the set of all
$\langle f, S,{\cal N},\Psi\rangle$ such that

(1) $S$ is a countable set of countable limit ordinals,

(2) $f$ is a function that
$S$-$*$-specializes $T$,

(3) for some
non-limit $\alpha=_{\rm def}{\rm lh}({\cal N})<\omega_1$ we have that
${\cal N}=\langle {\cal N}(i)\,\colon i<\alpha
\rangle$ is a tower (not necessarily continuous) of countable
elementary substructures of $H_\kappa$,

(4) for $i<j<{\rm lh}({\cal N})$ we have ${\cal N}(i)\in {\cal N}(j)$ and
if $\alpha\ne 0$ then $\{T,T^*,S^*,B\}\in{\cal N}(0)$, and
for every $i<{\rm lh}({\cal N})$ we have
$\omega_1\cap {\cal N}(i)\in S\}$,

(5) $\Psi$ is a
set of promises that $\langle f,S\rangle$ fulfills,

(6) for every $\beta<\omega_1$ we have that
$\{\Gamma\in\Psi\,\colon\allowbreak{\rm rk}({\rm min}(\Gamma))<\beta\}$
is countable,

(7) for every limit
ordinal $\alpha$ and $\gamma\in S^*$, if
$\{\omega_1\cap {\cal N}(\beta)\,\colon\allowbreak\beta<\alpha\}$ is
unbounded in $\gamma$ then $\alpha\in{\rm dom}({\cal N})$ and
${\cal N}(\alpha)=\bigcup\{{\cal N}(\beta\,\colon\allowbreak\beta<\alpha\}$,

(8) for all
$\beta\in{\rm dom}({\cal N})$, for all $x\in{\rm dom}(f)-{\cal N}(\beta)$ the following are
equivalent:

$i)$  $(\exists y<x)(y\in{\rm dom}(f)\cap{\cal N}(\beta)$ and $f(x)=f(y))$

$ii)$  $(\exists b\in B\cap{\cal N}(\beta))(x\in b)$

We order $P$ by declaring $\langle f, S,{\cal N},\Psi\rangle\leq
\langle f',S',{\cal N}',\Psi'\rangle$ iff
$S$ end-extends $S'$ and $f'\subseteq f$ and ${\cal N}$ end extends
${\cal N}'$ and
$\Psi'\subseteq\Psi$.

\proclaim Notations 9.5.  For $p\in P$ we give $f_p$, $S_p$, ${\cal N}_p$,
and $\Psi_p$  their obvious meanings, and
we set ${\rm ht}(p)= {\rm sup}(S_p)$ and $L_p=\bigcup\{{\cal N}_p(i)\,\colon
\allowbreak i\in{\rm dom}({\cal N}_p)\}$ and for $x\in L_p$ we let
$\rho_p(x)$ denote the least $\gamma$ such that
$x\in{\cal N}_p(\gamma)$. For $b\in B\cap L_p$ we set
$\delta_p(b)=\omega_1\cap{\cal N}_p(\rho_p(b))$ and we let
$\mu_p(b)$ denote the unique $x\in b$ such that
${\rm rk}(x)=\delta_p(b)$, and if $\delta_p(b)\in S_p$ then we 
 set $\sigma_p(b)$ equal to $f_p(\mu_p(b))$.  We set
 $U_p=\{x\in T\,\colon\allowbreak(\exists b\in B\cap L_p)\allowbreak
 (x\in b)\}$.

\proclaim Lemma 9.6.  Suppose $\lambda$ is a sufficiently large regular
cardinal (in particular, much larger than $\kappa$) and\/
 $M$ is a countable elementary substructure of
of\/ $H_\lambda$ containing $\{P,T,T^*,S^*,B,\kappa\}$.  Suppose
$p\in M\cap M$ and $n\in\omega$. Let $\delta=\omega_1\cap M$.
Suppose ${\overline x}\in T^n_\delta$ and suppose $R$ is
a finite rectangle and $n(R)=n$. 
 Suppose $z\in M\cap H_\kappa$. Then there is
$q\in P\cap M$ such that $q\leq p$ and $z\in L_q$ and 
$\heartsuit({\rm ht}(p),{\overline x}, f_q,R)$.

Proof: Let
 $N$ be a countable elementary substructure of $H_{\kappa}$ such that
 $N\in M$ and $\{z,R\cap
  \delta^n,T,T^*,S^*\}\in N$,
and let $\delta_\omega=\omega_1\cap N$.
Choose $\langle \delta_m\,\colon\allowbreak m<\omega\rangle$ an
  increasing sequence from $\delta_\omega\cap N$ cofinal in
  $\delta_\omega$ such that
  $\omega_1\cap L_p<\delta_0$ and
  $(\forall m\in\omega)\allowbreak(\delta_m\in
C(\Psi_p))$, where $C(\Psi_p)$ is as in the paragraph following Definition
7.7,  and such that
  for every $m\in\omega$, for every ${\overline y}\in\Gamma\cap
  T^{n(\Gamma)}_{\delta_m}$
  there is $W\subseteq\Gamma\cap T^{n(\Gamma)}_{\delta_{m+1}}$ such that
  for every $\{{\overline w},{\overline w}'\}\subseteq W$, if
  ${\overline w}\ne{\overline w}'$ then $\{{\overline w}(i)\,\colon
  \allowbreak i\in n(\Gamma)
  -G(\Gamma)\}$ is disjoint from $\{{\overline w}'(i)\,\colon\allowbreak
  i\in n(\Gamma)-G(\Gamma)\}$ and either $W$ is infinite or $G(\Gamma)=
  n(\Gamma)$.

  Let $\langle ({\overline y}_k,\Gamma_k,R'_k,t_k)\,
  \colon\allowbreak
  k\in\omega\rangle\in M$ list all quadruples $({\overline y},\Gamma,
R',t)$ such that
  $\Gamma\in\Psi_p$ and ${\overline y}\in\Gamma\cap N$ and
  $R'\subseteq\delta_\omega^{n(\Gamma)}$
  is a finite rectangle and $t\leq \omega$ and
  $(\forall i\in n(\Gamma))\allowbreak
  (i\in G(\Gamma)$ iff $(\exists b\in B\cap L_p)\allowbreak
  ({\rm min}(\Gamma)(i)\in b))$, with each such
  quadruple listed infinitely many times.

  Let $\langle x_i\,\colon\allowbreak i\in\omega\rangle\in M$ list
  $\bigcup\{T_{\delta_m}\,\colon\allowbreak m\leq\omega\}$.

  Working in $M$, build $\langle({\overline z}_m,Z_m,
  x^\#_m,X_m, f_m)\,\colon\allowbreak m\in\omega\rangle$ 
such that $f_0=f_p$ and for every $m\in\omega$ each of the following holds:

  (1) $f_m \subseteq f_{m+1}$

  (2) ${\rm dom}(f_{m+1})={\rm dom}(f_m)\cup\{x_m\}\cup Z_m\cup X_m$

  (3) if $\delta_{t_m}>{\rm rk}({\overline y}_m)$ then ${\overline y}_m<
{\overline z}_m$ and ${\rm rk}({\overline z}_m)=\delta_{t_m}$ and
$\{{\overline z}_m(i)\,\colon i\in G(\Gamma_m)\}= U_p\cap
{\rm range}({\overline z}_m)$ and
${\overline z}_m\in\Gamma_m$

  (4) $Z_m=\{{\overline z}_m(i)\,\colon\allowbreak i\in n(\Gamma_m)\}$

  (5) if $\delta_{t_m}\leq{\rm rk}({\overline y}_m)$ then $Z_m=\emptyset$

  (6)  if ${\rm rk}(x_m)\ne\delta_\omega$ or $x_m\in {\rm dom}(f_m)\cup Z_m$
  or there is no $b\in B\cap N$ such that $x_m\in b$
   then $X_m=\emptyset$; otherwise, $j_m\in\omega$
   is large
  enough that $(\forall x\in{\rm dom}(f_m)\cup Z_m)
  \allowbreak(x_m\restr\delta_{j_m}\not\leq x)$
   and
  $x^\#_m<x_m$ and ${\rm rk}(x^\#_m)=\delta_{j_m}$
  and $X_m=\{x^\#_m\}$

  (7) if $X_m\ne\emptyset$ then $f_{m+1}(x^\#_m)=f_{m+1}(x_m)$

  (8) for all $x\in{\rm dom}(f_{m+1})$, if there is
   $b\in B\cap L_p$ such that $x\in b$ then for the unique
  such $b$ we have $f_{m+1}(x) =\sigma_p(b)$

  (9) for every $x\in{\rm dom}(f_{m+1}-f_m)$, if
  there is no $b\in B\cap L_p$ such that $x\in b$ then
  $f_{m+1}(x)\notin\{f_{m+1}(x')\,\colon\allowbreak
   x'\in{\rm dom}(f_m)\cup\{x_m\}\cup Z_m\}$

  (10) for all $j\leq m$, if $\delta_{t_j}>{\rm rk}({\overline y}_j)$ then
  $\heartsuit({{\rm rk}({\overline y}_j,
  \overline z}_j, f_{m+1},R'_j)$

  (11)  $\heartsuit({\rm ht}(p),{\overline x},f_m,R)$

  There is no difficulty in meeting these requirements. 
  Set $q=\langle \bigcup\{f_m\,\colon\allowbreak m\in\omega\},\allowbreak
  S_p\cup\{\delta_m\,\colon m\leq\omega\},\allowbreak
  {\cal N}_p\hat{\  }\langle N\rangle,\allowbreak\Psi_p\rangle$.
  Then $q$ is as required in the conclusion of the Lemma.

\proclaim Lemma 9.7.  Suppose $T$ is an $\omega_1$-tree and $k\in\omega$ and
$\Delta\subseteq T^k$ is uncountable and downward closed, and suppose
${\overline x}\in T^k$ and every element of
$\Delta$ is comparable with ${\overline x}$.  Then there is
a promise $\Gamma\subseteq\Delta$ such that ${\rm min}(\Gamma) =
{\overline x}$.

Proof: We build $G\subseteq k$ and a sequence of uncountable
branches $\langle b_i\,\colon\allowbreak i\in G\rangle$ in stages as follows.

Initially, let $G_0=\emptyset$ and $\Delta_0=\Delta$
and $T^0=\{x\in T\,\colon\allowbreak
(\exists {\overline y}\in\Delta)\allowbreak(\exists i<k)\allowbreak
(x={\overline y}(i))\}$.

Stage $j$:

Case 1: $T^j$ is
 Aronszajn.

Take $G=G_j$. By Fact 7.5 we may take $\Gamma'\subseteq
\Delta_j$ such that for every $i\in k-G$ we have
${\rm min}(\Gamma)(i)={\overline x}(i)$.  Let $\Gamma=\{{\overline y}\in
\Delta\,\colon\allowbreak(\exists{\overline z}\in\Gamma')\allowbreak
((\forall m\in k-G)\allowbreak({\overline y}(m)={\overline z}(m)$ and
$(\forall m\in G)\allowbreak({\overline y}(m)$ is the unique element of
$b_m\cap T_{{\rm rk}({\overline y})}))\}$.

Case 2: Otherwise.

 Let $b$ be an uncountable branch of $T^j$ and fix
$i\in k-G_j$ such that $b\cap\{{\overline y}(i)\,
\colon\allowbreak{\overline y}\in\Delta_j
\}$ is uncountable.  Denote this $b$ by $b_i$.
Notice that because $\Delta_j$
 is downwards closed, we have
that $b_i$ is a subset of $\{{\overline y}(i)\,\colon\allowbreak
{\overline y}\in\Delta_j\}$.  Let $G_{j+1}=G_j\cup\{i\}$.
 For every ${\overline y}\in\Delta_j$ let
$s({\overline y})\in T^{k-G_{j+1}}$ be defined by
${\rm dom}(s({\overline y}))=k-G_{j+1}$ and for every $m\in k-G_{j+1}$
 we have
$s({\overline y})(m)={\overline y}(m)$. Set $\Delta_{j+1}=\{s({\overline y})
\,\colon
\allowbreak{\overline y}\in\Delta_j\}$, and let
$T^{j+1}=\{{\overline y}(m)\,\colon{\overline y}\in\Delta_{j+1}$ and
$m\in k-G_{j+1}\}$.  Now proceed to Stage $j+1$.

The Lemma is established.

\proclaim Lemma 9.8. Suppose $\lambda$ is a sufficiently large regular cardinal
and $M$ is a countable elementary substructure of $H_\lambda$ such that
$\{P,T,T^*,S^*,B,\kappa\}\in M$. Suppose $n\in\omega$ and
${\overline x}\in T^n_{\omega_1\cap M}$ and $R$ is a finite
rectangle with $n(R)=n$.  Suppose $p\in P\cap M$.
Then whenever $D\in M$ is a dense open subset of $P$,
 there is $q\leq p$ such that $q\in D\cap M$ and
$\heartsuit({\rm ht}(p),{\overline x},f_q,R)$.

Proof: Suppose $D$ is a counterexample.  We may assume $R\in M$
because if we replace each $R(i)$ with $R(i)\cap M$
we do not thereby change the truth of
$\heartsuit({\rm ht}(p),{\overline x},f,R)$ for any $f\in M$.
Set ${\overline z}={\overline x}\restr{\rm ht}(p)$ and set
$\Delta=\{{\overline y}\in T^n\,\colon\allowbreak
{\overline y}$ is comparable with ${\overline z}$ and there is no
$q\leq p$ such that $q\in D$ and ${\rm ht}(q)\leq{\rm rk}({\overline y})$
and $\heartsuit({\rm ht}(p),{\overline y},f_q,R)\}$.
Notice that $\Delta$ is downward closed and 
$\{{\overline y}\in T^n\,\colon\allowbreak{\overline y}<{\overline x}\}
\subseteq\Delta$. Necessarily $\Delta$ is uncountable because
$M\models``(\forall\alpha<\omega_1)\allowbreak(\Delta\cap T^n_\alpha\ne
\emptyset)$.''  Therefore by Lemma 9.7 we may take
$\Gamma\subseteq\Delta$ a promise with ${\rm min}(\Gamma)={\overline z}$.

Let $p'=\langle f_p,S_p,{\cal N}_p,\Psi_p\cup\{\Gamma\}\rangle$. Take
$r\leq p'$ such that $r\in D$.  Because
$\langle f_r,S_r,{\cal N}_r\rangle$ fulfills $\Gamma$,
we may take ${\overline w}\in\Gamma$ with
${\rm ht}(r)={\rm rk}({\overline w})$ and
$\heartsuit({\rm ht}(p),{\overline w},f_r,R)$.
Because ${\overline w}\in\Delta$, there is no $q\leq p$ such that $q\in D$ and
${\rm ht}(q)\leq{\rm rk}({\overline w})$ and $\heartsuit({\rm ht}(p),{\overline w},
f_q,R)$. But $r$ is a witness that there is such a $q$.  This
contradiction establishes the Lemma.

\proclaim Lemma 9.9. Suppose $\lambda$ is a sufficiently large regular cardinal
and $M$ is a countable elementary substructure of $H_\lambda$ such that
$\{P,T,T^*,S^*,B,\kappa\}\in M$. Suppose $n\in\omega$ and
${\overline x}\in T^n_{\omega_1\cap M}$ and $R$ is a finite
rectangle with $n(R)=n$.  Suppose $p\in P\cap M$.
Then 
 there is $q\leq p$ such that $q$ is $M$-generic and
 ${\rm ht}(q)=\omega_1\cap M$
 and
$\heartsuit({\rm ht}(p),{\overline x},f_q,R)$.

Proof: Take $N$ a countable elementary substructure of $H_\kappa$ such that
$N\in M$ and $\{p,\zeta,T,{\overline A}\cap\delta^n\}\in N$.
Set $\delta_\omega=\omega_1\cap N$ and choose $\langle\delta_m\,\colon\allowbreak
m\in\omega\rangle$ an increasing sequence from $\delta_\omega\cap N$ cofinal in
$\delta_\omega$, such that for every $\Gamma\in\Psi_p$ and every $m\in\omega$
 we have
$\delta_m\in C(\Gamma)$ and for every ${\overline y}\in\Gamma\cap T^{n(\Gamma)}_{\delta_m}$ there is
$W\subseteq\Gamma\cap T^{n(\Gamma)}_{\delta_{m+1}}$ such that for every
${\overline w}\in W$ we have ${\overline y}\leq{\overline w}$ and
for every $\{{\overline w},{\overline w}'\}\subseteq W$, if ${\overline w}\ne{\overline w}'$ then
$\{{\overline w}(i)\,\colon\allowbreak
i\in n(\Gamma)-G(\Gamma)\}$ is disjoint from
$\{{\overline w}'(i)\,\colon\allowbreak i\in n(\Gamma)-G(\Gamma)\}$, and either
$W$ is infinite or $G(\Gamma)=n(\Gamma)$.

Let $\langle\langle{\overline y}_k,\Gamma_k,{\overline A}^*_k,t_k\rangle\,\colon\allowbreak
k\in\omega\rangle$ list all $\langle {\overline y},\Gamma,{\overline A}^*,t\rangle$ such that $\Gamma\in\Psi_p$ and
${\overline y}\in\Gamma\cap N$ and ${\overline A}^*\subseteq(\delta_\omega)^n(\Gamma)$ is
a finite rectangle, and $t\leq\omega$, listed with infinitely many repetitions.

Let $S_q=S_p\cup\{\delta_t\,\colon\allowbreak t\leq\omega\}$.
Let $\langle x_m\,\colon\allowbreak m\in\omega\rangle$ list
$\bigcup\{T_{\delta_m}\,\colon\allowbreak m\leq\omega\}$.

Build $\langle f_m\,\colon\allowbreak m\in\omega\rangle$ such that $f_0=f_p$ and each of the following holds:

(1) $f_m\subseteq f_{m+1}$ and ${\rm dom}(f_{m+1})=
{\rm dom}(f_m)\cup\{x_m\}\cup Z_m\cup X_m$

(2) if $\delta_{t_m}>{\rm rk}({\overline y}_m)$ then ${\overline z}_m>{\overline y}_m$ and
${\rm rk}({\overline z}_m)=\delta_{t_m}$ and $\{{\overline z}_m(i)\,\colon\allowbreak i\in n(\Gamma_m)-G(\Gamma_m)\}$ is
disjoint from ${\rm dom}(f_m)$ and
${\overline z}_m\in\Gamma_m$ and
$Z_m=\{{\overline z}_m(i)\,\colon\allowbreak i\in n(\Gamma_m)\}$

(3) if $\delta_{t_m}\leq{\rm rk}({\overline y}_m)$ then $Z_m=\emptyset$

(4) if $t_m\ne\omega$ or $G(\Gamma_m)=\emptyset$
 then $X_m=\emptyset$; otherwise, $j_m\in\omega$ is
large enough that $x_m\restr\delta_{j_m}\not\leq x_j$ for all $j<m$, and
 $X_m=\{x^\#_m\}$ where $x^\#_m={\overline w}_m^\#(i)$ for some
$i\in G(\Gamma_m)$ and some
${\overline w}_m^\#\in\Gamma_m\cap T^{n(\Gamma_m)}_{\delta_{j_m}}$,
and ${\overline w}_m^\#(i)<x_m$ and
$(\forall b\in B\cap L_p)\allowbreak
({\overline w}_m^\#(i)\in b$ implies $x_m\in b)$

(5) if $X_m\ne\emptyset$ then $f_{m+1}(x^\#_m)=f_{m+1}(x_m)$

(6) if there is $b\in B\cap L_p$ such that $x_m\in b$ then for the unique such $b$ we have
$f_{m+1}(x_m)=\sigma_p(b)$

(7) if $\delta_m>{\rm rk}({\overline y}_m)$ then for all
$i<n(\Gamma_m)$, if for some $b\in B\cap L_p$ we have
${\overline z}_m(i)\in b$, then for the unique such $b$ we have
$f_{m+1}({\overline z}_m(i))=\sigma_p(b)$

(8) for every $j\leq m$ such that $\delta_{t_j}>{\rm rk}({\overline y}_j)$ we have
$\heartsuit({\overline z}_j,f_m,{\overline a}^*_j)$ implies
$\heartsuit({\overline z}_j,f_{m+1},{\overline a}^*_j)$,
and for all $i\in n(\Gamma_j)-G(\Gamma_j)$ we have
$f_{m+1}({\overline z}_j(i))\ne f_m(z)$ for all $z\in{\rm dom}(f_m)$ such that
$z$ is comparable with ${\overline z}_j(i)$

(9) if $x_m\notin{\rm dom}(f_m)\cup Z_m\cup X_m\cup(\bigcup(B\cap L_p))$ then
$f_{m+1}(x_m)\notin\{f_{m+1}(t)\,\colon\allowbreak
t\in{\rm dom}(f_m)\cup Z_m\cup X_m\}$

(10) $\heartsuit({\overline x},f_m,{\overline A})$ implies
$\heartsuit({\overline x},f_{m+1},{\overline A})$

Let $f_q=\bigcup\{f_m\,\colon\allowbreak m\in\omega\}$ and
${\cal N}_q={\cal N}_p\hat{ }\langle N\rangle$.

\section{Not adding reals}

In this section we discuss a sufficient condition for no reals
to be added.  This condition has two parts.
One part is a generalization of
[18, Definition 32], which is a variant of
Shelah's notion of ${\cal D}$-completeness
 [23, Chapter V]. The second part is
 Definition 32 given above.

\proclaim Lemma 41. Suppose  $\dot Q$ is
$(T,X)$-complete for $P$ and $\lambda$ is large for $\{P*\dot Q,X\}$ and
$M\prec N$ are countable elementary substructures of $H_\lambda$ and
$\{P*\dot Q, X, T\}\in M\in N$ and\/ $(p,\dot q)\in P*\dot Q\cap M$
and $G\in\Gen(M,T,P,p)\cap N$.
Then there are a $P$-name $\dot s$ and a set $G'\in
\Gen(M,T,P*\dot Q,(p,\dot q))$ such that
$G=\{p'\in P\,\colon\allowbreak(\exists \dot r)((p',\dot r)\in G')\}$
and whenever $\tilde p$ is a lower bound for $G$
which is  $(M,P,T)$-completely preserving
and $(N,P,T)$-preserving then
$(\tilde p,\dot s)$ is an $(M,P,T)$-completely
preserving lower bound for $G'$.

Proof:  Let $G'$ be as in the conclusion of Definition 39 and for
every $p'\in P$ such that $p'$ is a lower bound for
$G$ which is both $(M,P,T)$-completely preserving and
$(N,P,T)$-preserving
 let $\dot s(p')$ be as in the conclusion of Definition 39.
Let $\calj$ be a maximal antichain of the set of such $p'$ and take
$\dot s$ such that $(\forall p'\in\calj)\allowbreak(p'\forces``\dot s=
\dot s(p')$''). We have that $G'$ and $\dot s$ are as required.

\proclaim Definition 42.  Suppose $\langle P_\eta\,\colon
\allowbreak\eta\leq\alpha\rangle $ is a countable support iteration
and\/ $T$ is Suslin and $X$ is any set. We say that
$P_\alpha$ is\/ {\rm $(T,X)$-strictly complete} iff\/ {\bf
whenever} $\lambda$ is large for $\{P_\alpha,X\}$ and
$M$ is a countable elementary substructure of $H_\lambda$ and
$\{P_\alpha,X,T\}\in M$ and $\alpha^*$ is the order-type of
$\alpha\cap M$ and 
${\cal N}=\langle N_i\,\colon\allowbreak
i\leq\alpha^*\rangle$ is a $\lambda$-tower
for $M$ and $p\in P_\alpha\cap M$ and
 $\eta\in\alpha\cap M$ and $\eta^*$ is the order-type of\/
 $\eta\cap M$ 
and $G\in\Gen(M,P_\eta,p\restr\eta)\cap N_{\eta^*+1}$
{\bf then} there are $G'\in\Gen(M,P_\alpha,p)$ and a $P_\eta$-name
$\dot s$ 
such that $\{r\restr\eta\,\colon\allowbreak r\in G'\}=G$ and
{\tt whenever} $\tilde p$ is an $(M,P_\eta,T)$-completely preserving
 lower bound for $G$ and
$\tilde p$ is
$(\langle
 N_i\,\colon\allowbreak\eta^*<i\leq\alpha^*\rangle,
 \allowbreak P_\eta,T)$-preserving 
{\tt then} we have that
there is $\tilde s\in P_\alpha$ such that
$\tilde s\restr\eta=\tilde p$ and
$\tilde s$ is an $(M,P_\alpha,T)$-completely
preserving lower bound for $G'$
and $\tilde p\forces``\tilde s\restr[\eta,\alpha)=\dot s$''
and $\supt(\tilde s)\subseteq\eta\cup N_{\alpha^*}$.

\proclaim Lemma 43. Suppose $P_\alpha$ is $(T,X)$-strictly complete
for some $X$.
Then $P_\alpha$ does not add reals.

Proof: Simply take $\eta=0$ in Definition 42.

\proclaim Lemma 44. Suppose
 $\langle P_\eta\,\colon\allowbreak\eta\leq\alpha
\rangle$ is a countable support iteration
based on $\langle\dot Q_\eta\,\colon\allowbreak
\eta<\alpha\rangle$ and $T$ is Suslin
 and  for all $\eta<\alpha$ we have
$\dot Q_\eta$ is $(T,X_\eta)$-complete for $P_\eta$,
and suppose for every $\beta\leq\alpha$
we have that $P_\beta$ is strictly strongly $T$-preserving.
Then $P_\alpha$ is $(T,X)$-strictly complete where
$X=\langle X_\eta\,\colon\allowbreak \eta<\alpha\rangle$.

Proof:
We work by induction on $\alpha$.
Let $\lambda$, $M$, $\alpha^*$, $\eta$, $\eta^*$,
${\cal N}=
\langle N_i\,\colon\allowbreak i\leq\alpha^*\rangle$,
$p$, and $G$ be as in the hypothesis of Definition 42.

Suppose first that $\alpha=\beta+1$.
Let $\beta^*$ be the order-type of $\beta\cap M$.
By the induction hypothesis we may take
$G_1\in\Gen(M,P_{\beta},p\restr\beta)$ and $\dot s_1$
such that $\{r\restr\eta\,\colon\allowbreak r\in G_1\}=G$ and
whenever $\tilde p$ is a lower bound for $G$ which is both
$(M,P_\eta,T)$-completely preserving and
$(\langle N_i\,\colon\allowbreak \eta^*<
i\leq\beta^*\rangle,\allowbreak P_\eta,T)$-preserving
then we have that there is $ s^*\in P_{\beta}$  such that
$s^*$ is an $(M,P_\beta,T)$-completely preserving lower bound for $G_1$ and
$s^*\restr\eta=\tilde p$ and
$\tilde p\forces``s^*\restr[\eta,\beta)=\dot s_1$''
and $\supt(s^*)\subseteq\eta\cup N_{\beta^*}$.
By elementarity, we may assume that $G_1$ and $\dot s_1$ are
elements of $N_{\alpha^*}$.
Because $P_\beta$ is strictly strongly
$T$-preserving, we may take $\dot s'_1$ such that whenever
$\tilde p$ is a lower bound for $G$ which is 
both $(M,P_\eta,T)$-completely
preserving and $(\langle N_i\,\colon\allowbreak
\eta^*<i\leq\alpha^*\rangle,\allowbreak
P_\eta,T)$-preserving  then
$\tilde p\forces``\dot s'_1\leq\dot s_1$'' and $(\tilde p,\dot s_1)$ is
$(N_{\alpha^*},P_\beta,T)$-preserving and
$\tilde p\forces``\supt(\dot s'_1)\subseteq N_{\alpha^*}[G_{P_\eta}]$.''
Necessarily we have that $(\tilde p,\dot s_1)$ is
$(M,P_\beta,T)$-completely preserving.
By Lemma 41 we may take $\dot s_2$ and
$G'\in\Gen(M,P_\alpha,p)$ such that
$G_1=\{r\restr\beta\,\colon\allowbreak r\in G'\}$ and whenever
$\tilde p$ is a lower bound for $G_1$
which is both $(M,P_\beta,T)$-completely preserving and
 $(N_{\alpha^*},P_\beta,T)$-preserving then $(\tilde p,\dot s_2)$ is an
 $(M,P_\alpha,T)$-completely preserving
lower bound for $G'$. Let $\dot s$ be the $P_\eta$-name
for the pair $(\dot s_1',\dot s_2)$. Then
$\dot s$ and $G'$ are as required.

Now we consider the case where $\alpha$ is a limit ordinal.
Let $\langle\alpha_n\,\colon\allowbreak n\in\omega\rangle$
be an increasing sequence from $\alpha\cap M$ cofinal in
$\sup(\alpha\cap M)$ such that $\alpha_0=\eta$.
For every integer $n\geq 0$ let
$\alpha^*_n$ be the order-type of $\alpha_n\cap M$.
 Let $\langle\tau_n\,\colon\allowbreak n\in\omega\rangle$
list the set of all $P_\alpha$-names $\tau$ in $M$ such that
$\bfone\forces``\tau$ is an ordinal.''
Let $\langle\langle x_n,A_n\rangle\,\colon\allowbreak n\in\omega\rangle$
list the set of all pairs $\langle x,A\rangle$ such that
$x\in T$ and $\rk(x)=\omega_1\cap M$ and
$A\in M$ and $A$ is a $P_\alpha$-name for a subset of $T$.

Build $\langle\langle G_n,\dot s_n,\dot s'_n,p_n\rangle\,\colon\allowbreak
n\in\omega\rangle$ such that $G_0=G$ and $p_0=p$ and
each of the following:

\noindent(1) $G_n\in\Gen(M,P_{\alpha_n},p_n\restr\alpha_n)\cap
N_{\alpha^*_n+1}$

\noindent(2) $p_{n+1}\leq p_n$ and $p_{n+1}\in P_\alpha\cap M$ and
$p_{n+1}\restr\alpha_n\in G_n$ and
$p_{n+1}\restr\alpha_n\forces``p_{n+1}\restr[\alpha_n,\alpha)$ decides
the value of $\tau_n$, and either $x_n\notin A_n$ or there are
$y<x_n$ and $z\in T\cap M$ such that
$z\not<x_n$ and $p_{n+1}\restr[\alpha_n,\alpha)\forces`
\{y,z\}\subseteq A_n$.'\thinspace''

\noindent(3) $\dot s_n\in N_{\alpha^*_{n+1}+1}$ and
whenever $\tilde p$ is a lower bound for $G_n$ and $\tilde p$ is
both $(M,P_{\alpha_n},T)$-completely preserving and
$(\langle N_i\,\colon\allowbreak\alpha^*_n<i\leq\alpha^*_{n+1}\rangle,
\allowbreak P_{\alpha_n},T)$-preserving,
then there is $\tilde s\in P_{\alpha_{n+1}}$ such that
$\tilde s\restr\alpha_n=\tilde p$ and $\tilde p\forces``\tilde s
\restr[\alpha_n,\alpha_{n+1})=\dot s_n$'' and
$\tilde s$ is an $(M,P_{\alpha_{n+1}},T)$-competely preserving
 lower bound for $G_{n+1}$.

\noindent(4) $G_n=\{r\restr\alpha_n\,\colon\allowbreak r\in G_{n+1}\}$.

\noindent(5) whenever $\tilde p$ is a lower bound for $G_n$ and
$\tilde p$ is both $(M,P_{\alpha_n},T)$-completely preserving and
$(\langle N_i\,\colon\allowbreak
\alpha^*_n<i\leq\alpha^*\rangle,\allowbreak
P_{\alpha_n},T)$-preserving, then $\tilde p\forces
``\dot s'_n\leq\dot s_n$'' and there
 is $\tilde s\in P_{\alpha_{n+1}}$ such that
$\tilde s\restr\alpha_n=\tilde p$ and $\tilde p\forces``\tilde s
\restr[\alpha_n,\alpha_{n+1})=\dot s'_n$'' and
$\supt(\tilde s)\subseteq\alpha_n\cup N_{\alpha^*}$ and
$\tilde s$ is
$(\langle N_i\,\colon\allowbreak\alpha^*_{n+1}<
i\leq\alpha^*\rangle,\allowbreak
P_{\alpha_{n+1}},T)$-preserving (necessarily, $\tilde s$ is
$(M,P_{\alpha_{n+1}},T)$-completely preserving).

The construction proceeds as follows. Given $G_n$ and $p_n$,
construct $p_{n+1}$ as in (2) as follows. Choose $\dot q\in M$ such that
$p_n\restr\alpha_n\forces``\dot q\leq p_n\restr[\alpha_n,\alpha)$
and $\dot q$ decides the value of $\tau_n$.''
Let $E=\{r\leq p_n\restr\alpha_n\,\colon
\allowbreak(\exists s\in P_\alpha)\allowbreak
(s\restr\alpha_n=r$ and $r\forces``s\restr[\alpha_n,\alpha)=\dot q$''$)\}$.
Because $E\in M$ we may take $r_1\in E\cap G_n$. Take
$q_1\in P_\alpha\cap M$ such that $q_1\restr\alpha_n=r_1$ and
$r_1\forces``q_1\restr[\alpha_n,\alpha)=\dot q$.'' Let
$X=\{w\in T\,\colon\allowbreak q_1\notforces``w\notin A_n$''$\}$.
We select $q_3$ as follows. 
If $x_n\notin X$ let $q_3=q_1$. Otherwise, take $y<x_n$ such that
$y\in X$. We have $r_1\forces``
q_1\restr[\alpha_n,\alpha)\notforces`y\notin A_n$'\thinspace''
so we may take $\dot q_2\in M$ such that
$r_1\forces``\dot q_2\leq q_1\restr[\alpha_n,\alpha)$ and
$\dot q_2\forces`y\in A_n$.'\thinspace''
Let $E_1=\{r\leq r_1\,\colon\allowbreak(\exists s\in P_\alpha)\allowbreak
(s\restr\alpha_n=r$ and $r\forces``s\restr[\alpha_n,\alpha)=\dot q_2$''$)\}$.
 Because $E_1\in M$ we may take
$r_2\in E_1\cap G_n$. Then take $q_3\in P_\alpha\cap M$
such that $q_3\restr\alpha_n=r_2$ and
$r_2\forces``q_3\restr[\alpha_n,\alpha)=\dot q_2$.''
Let $Y=\{w\in T\,\colon\allowbreak
q_3\notforces``w\notin  A_n$''$\}$. We build $p_{n+1}$ as follows.
If $x_n\notin Y$ then we let $p_{n+1}=q_3$. Otherwise,
 take $z\in Y\cap M$ such that
$z\not<x_n$.
 We have $r_2\forces``q_3\restr[\alpha_n,\alpha)
 \notforces`z\notin A_n$'\thinspace''
so we may take $\dot q_4\in M$ such that
$r_2\forces``\dot q_4\leq q_3\restr[\alpha_n,\alpha)$ and
$\dot q_4\forces`z\in A_n$.'\thinspace''
Let $E_2=\{r\leq r_2\,\colon\allowbreak(\exists s\in P_\alpha)\allowbreak
(s\restr\alpha_n=r$ and $r\forces``s\restr[\alpha_n,\alpha)=\dot q_4$''$)\}$.
 Because $E_2\in M$ we may take
$r_3\in E_2\cap G_n$. Then take $p_{n+1}\in P_\alpha\cap M$
such that $p_{n+1}\restr\alpha_n=r_3$ and $r_3\forces``
p_{n+1}\restr[\alpha_n,\alpha)=\dot q_4$.''

Given $p_{n+1}$, use the fact that $P_{\alpha_{n+1}}$ is
$(T,X)$-strictly complete to take $G_{n+1}$ and
$\dot s_n$ as in (1) and (3) and (4).
Finally, use Lemma 36 to take $\dot s'_n$ as in (5).

Let $G'=\{p\in M\,\colon\allowbreak(\exists n\in\omega)\allowbreak
(p_n\leq p)\}$, and let
$\dot s$ be the $P_\eta$-name for the concatenation of
$(\dot s'_0,\dot s'_1,\ldots)$, followed by $\bfone_{\zeta,\alpha}$
where $\zeta=\sup(\alpha\cap M)$.

We show that this choice of $G'$ and $\dot s$ works.
Given $\tilde p$ a lower bound for $G$ which is
both $(M,P_\eta,T)$-completely
preserving and $(\langle N_i\,\colon\allowbreak
\eta^*<i\leq\alpha^*\rangle,\allowbreak
P_\eta,T)$-preserving,
we build $\langle \tilde p_n\,\colon\allowbreak
n\in\omega\rangle$ such that $\tilde p_0=\tilde p$ and
for every $n\in\omega$ we have
$\tilde p_{n+1}\restr\alpha_n=\tilde p_n$ and
$\tilde p_{n+1}$ is a lower bound for $G_{n+1}$ and
$\tilde p_{n+1}$ is both
$(M,P_{\alpha_{n+1}},T)$-completely
preserving and $(\langle N_i\,\colon\allowbreak
\alpha^*_{n+1}<i\leq\alpha^*\rangle,\allowbreak
P_{\alpha_{n+1}},T)$-preserving, and
$\tilde p_n\forces``\tilde p_{n+1}\restr[\alpha_n,\alpha_{n+1})=
\dot s'_n$'' and $\supt(\tilde p_{n+1})\subseteq\eta\cup N_{\alpha^*}$.
This is possible because given $\tilde p_n$, there is a
$P_{\alpha_n}$-name $E$ such that
$\tilde p_n\forces``E$ is a closed subset of
$\Spec(\langle N_i\,\colon\allowbreak\alpha^*_n<i\leq\alpha^*\rangle,
\allowbreak P_{\alpha_n},T)$ of order-type
$(\alpha^*+1)-(\alpha^*_n+1)$ and
$(\forall i\in E)\allowbreak(\forall j\in i\cap E)\allowbreak
(E\cap N_j\in N_i[G_{P_{\alpha_n}}])$.''
Because $\tilde p_n\forces``(\alpha^*_{n+1}+1)\cap E$ has
order-type at most $(\alpha^*_{n+1}+1)-(\eta^*+1)$,''
we have $\tilde p_n\forces``\{i\in E\,\colon\allowbreak
\alpha^*_{n+1}<i\}$ has order-type at least
$(\alpha^*+1)-(\alpha^*_{n+1}+1)$, and hence has order-type
exactly equal to $(\alpha^*+1)-(\alpha^*_{n+1}+1)$.''
Hence we may proceed to take $\tilde p_{n+1}$ as given
above.

Let $r\in P_\alpha$ be such that
$\supt(r)\subseteq\sup(\alpha\cap M)$ and
$(\forall n\in\omega)\allowbreak( r\restr\alpha_n=\tilde p_n)$.
We have $\tilde p\forces``r\restr[\eta,\alpha)=\dot s$'' and
$r$ is an $(M,P_\alpha,T)$-completely preserving lower bound for $G'$.

\end{document}